\documentclass[10pt,twoside]{article}
\usepackage{amsmath,amsthm}
\usepackage[latin1]{inputenc}
\usepackage{picinpar}
\usepackage{longtable}
\usepackage{amsrefs}
\usepackage[dvips]{graphicx}
\usepackage[mathscr]{eucal}
\usepackage{calc}
\usepackage{ifthen}
\usepackage{dcpic,pictexwd}
\usepackage[all]{xy}
\usepackage{enumerate}
\usepackage{amssymb,latexsym}
\usepackage{amsthm}
\usepackage[dvips]{graphics}
\usepackage{color}
\usepackage{tabularx}
\usepackage{mathtools}
\usepackage{makeidx}

\usepackage{textcomp}

\theoremstyle{definition}
\newtheorem{teor}{Theorem}
\newtheorem{lema}{Lemma}
\newtheorem{prop}{Proposition}
\newtheorem{corol}{Corollary}

\newtheorem{Nota}{Remark}

\begin{document}
\thispagestyle{plain}
\par\bigskip

\begin{center}
\textbf{Coxeter's Frieze Patterns Arising from Dyck Paths}
\end{center}

\par\bigskip

\begin{centering}
Agust\'{\i}n Moreno Ca\~{n}adas\\
Isa\'{\i}as David Mar\'{\i}n  Gaviria\\
Gabriel Bravo Rios \\
Pedro Fernando Fern\'{a}ndez Espinosa\\
\end{centering}

\par\bigskip

\par\bigskip

\begin{centering}

\textbf{Abstract}\par\bigskip

\end{centering}

Frieze patterns are defined by objects of a category of Dyck paths, to do that, it is introduced the notion of diamond of Dynkin type $\mathbb{A}_{n}$. Such diamonds constitute a tool to build integral frieze patterns.\par\bigskip

\small{}

\par\bigskip



\section{Introduction}
Frieze patterns (as shown below) were introduced by Coxeter  in the early 70s \cite{Coxeter}. According to Propp \cite{Propp}, they arose from Coxeter's study of metric properties of polytopes, and served as useful scaffolding for various sorts of metric data.\\

\[\xymatrix@=1.5mm{\dots && 0 && 0 && 0 && 0 &&\dots\\
& 1 &&1 && 1 && 1 && 1  &\\
\dots && m_{-1,-1} && m_{00} && m_{11} && m_{22} &&\dots\\
& m_{-2,-1} &&m_{-1,0} && m_{0,1} && m_{12} && m_{23} &\\
\dots && \dots && \dots && \dots && \dots  &&\dots\\
& 1 &&1 && 1 && 1 && 1 &\\
\dots && 0 && 0 && 0 && 0 && \dots\\
}\]\\

Such patterns are defined as grids of numbers bounded from above by an infinite row of 0s followed by a row of 1s and such that every four adjacent numbers of the following form
\[\xymatrix@=2mm{ & b&\\
a&&c\\
&d&\\
}\]
satisfy the arithmetic (or frieze) rule 
\begin{equation}\label{coxetereq}
\begin{split}
ac-bd=1.
\end{split}
\end{equation}

The third line (i.e., the first nontrivial row) can be chosen in an arbitrary way and then complete using the rule (\ref{coxetereq}), which was named the modular equation by Coxeter \cite{Coxeter}.\\

A frieze is called closed, if it is also bounded from below by a line of 1s (followed by a line of 0s). A frieze is called integral if it consists of positive integers. 
The sequence of integers  in the first non-trivial row  is called \textit{quiddity sequence}. This sequence completely determines the frieze pattern. Each frieze pattern is also periodic since it is invariant under glide reflection. The order of the frieze pattern is defined to be the number of rows minus one. It follows that each frieze pattern of order $n$ is $n$-periodic \cite{Baur1}. Conway and Coxeter  classified completely the frieze patterns whose entries are positive integers, and show that these frieze patterns constitute a manifestation of the Catalan numbers \cite{Conway, Conway1} giving a bijection between positive integer frieze patterns and triangulations of regular polygons with labeled vertices. \\

Frieze patterns appear independently in the 70s in the context of quiver representations, in such a case, the local arithmetic rule is an additive analogue of Coxeter\textquotesingle s unimodular rule. The generalization of the Coxeter\textquotesingle s unimodular rule on Auslander-Reiten quivers was found by Caldero and Chapoton \cite{Caldero1}. Assem,  Reutenauer and  Smith introduced also a generalization of friezes by associating frieze patterns to Cartan matrices \cite{Assem1}.\\

According to Morier-Genoud \cite{Morier} there are mainly three approaches for the study of friezes:

\begin{enumerate}
\item A representation theoretical and categorical approach, in deep connection with the theory of cluster algebras, where entries in the friezes are rational functions.
\item A geometrical approach, in connection with moduli spaces of points in projective space and Grassmannians, where entries in the friezes are more often real or complex numbers.
\item A combinatorial approach, focusing on friezes with positive integer entries.

\end{enumerate}

Many authors have studied friezes from the different points of view finding connections with different branches of mathematics \cite{Morier}. For instance, Baur et al \cite{Baur1} studied mutation of friezes proving how mutation of a cluster affects the associated frieze. On the other hand, Fontaine and Plamondon \cite{Fontaine} presented a formula for the number of friezes of type $\mathbb{B}_{n}, \mathbb{C}_{n}, \mathbb{D}_{n}$, and $\mathbb{G}_{2}$. They conjectured that the number of friezes of type $\mathbb{E}_{6}, \mathbb{E}_{7}, \mathbb{E}_{8}$ and $\mathbb{F}_{4}$ is 868, 4400, 26592 and 112, respectively. In this way, the number of friezes can be defined as a Dynkin function in the sense of Ringel \cite{Ringel}.\\

In this paper, frieze patterns are interpreted as objects of a novel category of Dyck paths introduced recently by Ca\~{n}adas and Rios \cite{Canadas}.\\

The following is a list of our main results, all of them dealing with integral closed friezes.\\

\begin{enumerate}
\item It is introduced the notion of diamond of Dynkin type $\mathbb{A}_{n}$, and some of its properties are proved (see, Proposition \ref{prop4.1}, Proposition \ref{prop4.2} and Theorem \ref{teor4.1}). Such diamonds are used as a tool to build frieze patterns.

\item It is proved that there is a bijective correspondence between the set of all vectors associated to positive integral diamonds of Dynkin type $\mathbb{A}_{n}$ and triangulations of a polygon with $n+3$ vertices. This result is a consequence of a bijection between such triangulations and Dyck paths of length $2(n+1)$ (see Lemma \ref{lema4.3} and Theorem \ref{teor4.2}).

\item It is proved that if  $\mathscr{C}(A^{0},t)=\lbrace A^{t} \rbrace_{0 \leq t \leq p-1}$ is the minimal $p$-cycle  generated by a diamond $A^{0}$ of Dynkin type $\mathbb{A}_{n}$. Then $\mathscr{C}(A^{0},t)$  is in  surjective correspondence with a direct sum of $p$ indecomposable objects of a Dyck paths category (see Theorem \ref{main}).

\end{enumerate}

The paper is distributed as follows. In section 2, we recall main definitions and notation to be used throughout the paper, in particular, it is reminded the notion of a Dyck paths category. In section 3, we present the main results, where it is introduced the notion of diamond of Dynkin type $\mathbb{A}_{n}$ and some of its properties.


\section{Preliminaries}

In this section, we recall main definitions and notation to be used throughout the paper~\cites{Barcucci, Canadas, Gunawan, Stanley}.

\subsection{Cluster Algebras From Quivers}

 For quivers, cluster algebras are defined as follows:\par\bigskip

 Fix an integer $n\geq1$. In this case, a seed $(Q, u)$ consists of a finite quiver $Q$ without loops or 2-cycles with vertex set $\{1,\dots, n\}$, whereas $u$ is a free-generating  set $\{u_{1},\dots, u_{n}\}$ of the field $\mathbb{Q}(x_{1},\dots, x_{n})$.\par\bigskip

 Let $(Q, u)$ be a seed and $k$ a vertex of $Q$. The mutation $\mu_{k}(Q, u)$ of $(Q, u)$ at $k$ is the seed $(Q',u')$, where;

 \begin{enumerate}[(a)]
 
\item $Q'$ is obtained from $Q$ as follows;

\begin{enumerate}[(1)]
\item reverse all arrows incident with $k$,
\item for all vertices $i\neq j$ distinct from $k$, modify the number of arrows between $i$ and $j$, in such a way that a system of arrows of the form $(i \stackrel{r}{\longrightarrow} j, i \stackrel{s}{\longrightarrow} k, k \stackrel{t}{\longrightarrow} j) $ is transformed into the system $(i \stackrel{r+st}{\longrightarrow} j, k \stackrel{s}{\longrightarrow} i, j \stackrel{t}{\longrightarrow} k) $. And the system $(i \stackrel{r}{\longrightarrow} j, j\stackrel{t}{\longrightarrow} k, k \stackrel{s}{\longrightarrow} i) $ is transformed into the system $(i \stackrel{r-st}{\longrightarrow} j, i \stackrel{s}{\longrightarrow} k, k \stackrel{t}{\longrightarrow} j) $.  Where, $r$, $s$ and $t$ are non-negative integers, an arrow $i \stackrel{l}{\longrightarrow} j$, with $l\geq0$ means that $l$ arrows go form $i$ to $j$ and an arrow $i \stackrel{l}{\longrightarrow} j$, with $l\leq0$ means that $-l$ arrows go from $j$ to $i$.
\end{enumerate} 

\item $u'$ is obtained form $u$ by replacing the element $u_{k}$ with

\begin{equation}\label{exchangeruleq}
 \begin{split}
 u_{k}&=\frac{1}{u_{k}}\underset{\mathrm{arrows}\hspace{0.1cm}i\rightarrow k}{\prod}u_{i}+\underset{\mathrm{arrows}\hspace{0.1cm}k\rightarrow j}{\prod}u_{j}.
 \end{split}
 \end{equation}

 \end{enumerate}

 If there are no arrows from $i$ with target $k$, the product is taken over the empty set and equals 1. It is not hard to see that $\mu_{k}(\mu_{k}(Q, u))=(Q, u)$.
 In this case the matrix mutation $B'$ has the form

  \[b'_{ij}=
\begin{cases}
-b_{ij},  &   \hspace{0.2cm}\text{if}\hspace{0.1cm}i=k\hspace{0.1cm}\text{or}\hspace{0.1cm} j=k, \\
b_{ij}+sgn(b_{ik})[b_{ik}b_{k j}]_{+},  &   \hspace{0.2cm}\text{else}, 
\end{cases}\]

where $[x]_{+}=\mathrm{max(x,0)}$. Thus, if $Q$ is a finite quiver without loops or 2-cycles with vertex set $\{1,\dots, n\}$, the following interpretations take place:

\begin{enumerate}
\item the clusters with respect to $Q$ are the sets $u$ appearing in seeds, $(Q, u)$ obtained from a initial seed $(Q, x)$ by iterated mutation,
\item the cluster variables for $Q$ are the elements of all clusters,
\item the cluster algebra $\mathscr{A}(Q)$ is the $\mathbb{Q}$-subalgebra of the field $\mathbb{Q}(x_{1},\dots, x_{n})$ generated by all the cluster variables.
\end{enumerate}

 As example, the cluster variables associated to the quiver $Q=1\longrightarrow2$ are:\par\bigskip 

 \begin{centering}
 $\{x_{1}, x_{2}, \frac{1+x_{2}}{x_{1}}, \frac{1+x_{1}+x_{2}}{x_{1}x_{2}}, \frac{1+x_{1}}{x_{2}}\}$.\par\bigskip
 \end{centering}

Regarding cluster algebras arising from quivers, we recall that, Fomin and Zelevinsky \cite{Fomin1} proved that any cluster algebra $\mathcal{A}(Q)$ of finite type has a finite set of cluster variables and the following result.

\begin{teor}
\textit{The cluster algebra $\mathcal{A}(Q)$ is of finite type if and only if $Q$ is mutation-equivalent to an orientation of a simply-laced Dynkin diagram, $\mathbb{A}_{n}, n\geq 1$, $\mathbb{D}_{n}, n\geq 4$, $\mathbb{E}_{6}$, $\mathbb{E}_{7}$ and $\mathbb{E}_{8}$.}
\end{teor}

\subsection{Friezes}

An alternative way to define friezes is to say that they are ring homomorphisms from a cluster algebra to the ring of integers such that all cluster variables are sent to positive integers  \cite{Fontaine}. Let $Q$ be a quiver without loops and $2$-cycles and let $\mathcal{A}(Q)$ be the corresponding cluster algebra  with trivial coefficients \cite{Gunawan}, then: 

\begin{enumerate}
\item [(i)] A \textit{frieze of type} $Q$ is a ring homomorphism $\mathcal{F}: \mathcal{A}(Q)\rightarrow \mathbf{R}$ from the cluster algebra to an integral domain $\mathbf{R}$. The frieze is called \textit{integral} if $\mathbf{R}=\mathbb{Z}$.
 \item [(ii)] An integral frieze is said to be \textit{positive} if every cluster variable in $\mathcal{A}(Q)$ is mapped by $\mathcal{F}$ to a positive integer.
 \end{enumerate}

Let \textbf{x}$=(x_{1}, \dots, x_{n})$ be a cluster of $\mathcal{A}(Q)$, then:

\begin{enumerate}
\item [(iii)] A vector $(a_{1}, \dots, a_{n}) \in \mathbf{R}^{n}$ is called a \textit{frieze vector} relative to \textbf{x} if the frieze $\mathcal{F}$ defined by the assignment $\mathcal{F}(x_{i})=a_{i}$ has values in $\mathbf{R}$, $\mathcal{F}$ is said to be \textit{unitary} if there exists a cluster $\textbf{x}$ such that $\mathcal{F}(x)$ is a unit in $\textbf{R}$, for all $x\in\textbf{x}$. If the frieze $\mathcal{F}$ is unitary we say that the frieze vector $(a_{1}, \dots, a_{n})$ is \textit{unitary}.
\item [(iv)]A vector $(a_{1}, \dots, a_{n}) \in \mathbb{Z}_{>0}^{n}$ is called a \textit{positive} frieze vector relative to \textbf{x} if the frieze $\mathcal{F}$ defined by $\mathcal{F}(x_{i})=a_{i}$ is positive integral.
\end{enumerate}
The following is an example of a frieze pattern. Hereinafter, frieze patterns are assumed to be integral closed friezes.\\

\[\xymatrix@=2mm{\dots && 0 && 0 && 0 && 0 && 0 &&\dots\\
& \dots &&1 && 1 && 1 && 1 && \dots &\\
\dots && 2 && 3 && 1 && 2 && 3 &&\dots\\
& \dots &&5 && 2 && 1 && 5 && \dots &\\
\dots && 2 && 3 && 1 && 2 && 3 &&\dots\\
& \dots &&1 && 1 && 1 && 1 && \dots &\\
\dots && 0 && 0 && 0 && 0 && 0 &&\dots\\
}\]\\

\subsection{Dyck Paths Categories}

In this section we recall the definition and main properties of the category of Dyck paths as Ca\~nadas and Rios describe in \cite{Canadas}.\\

A Dyck path is a lattice path in $\mathbb{Z}^2$  from $(0,0)$ to $(n, n)$ with steps $(1,0)$ and $(0,1)$ such that the path never passes below the line $y=x$. The number of Dyck paths of length $2n$ is equal to $C_{n}=\frac{1}{n+1}\binom{2n}{n}$,  the nth Catalan number  \cite{Stanley}.\\

The set of \textit{Dyck words} is the set of words $w$ in the free monoid $X^{\ast}=\lbrace U,D\rbrace^{\ast}$ satisfying the following two conditions \cite{Barcucci}:

\begin{itemize}
\item for any left factor $u$ of $w$ ( i.e., $w=uv$ for some suitable word $v$), $|u|_{U}\geq |u|_{D}$,
\item $|w|_{U}=|w|_{D}$,
\end{itemize}

where $|w|_{a}$ is the number of occurrences of the letter $a \in X=\lbrace U,D\rbrace$ in the word $w$.\\

Henceforth, Dyck words defined as before are used to denote Dyck paths.\\

Let $\mathfrak{D}_{2n}$ be the set of all Dyck paths of length $2n$, let  $UWD=Uw_{1}\dots w_{n-1}D$ be a Dyck path in $\mathfrak{D}_{2n}$ with $A=\lbrace UD, DU, UU, DD\rbrace$ being  the set of all possible choices in $W$.\\

The \textit{support} of $UWD$ (denoted by $\text{Supp } UWD\subseteq\{1,2,\dots, n-1\}=\textbf{n-1}$) is a set of indices (of the $w_{i}$s) such that \[\text{Supp}\hspace{0.1cm}UWD= \{ q\in\textbf{n-1} \text{ } | \text{ } w_{q}=UD \text{ or }w_{q}=UU  \text {, } 1 \leq q \leq n-1 \}. \]

A map $f:A\longrightarrow A$ such that for any $w\in A$, it holds that $f(w)= f(ab)=w^{-1}=ba$, $a, b\in\{U, D\}$ is said to be a \textit{shift}. An \textit{unitary shift} is a map $f_{i}:\mathfrak{D}_{2n}\longrightarrow \mathfrak{D}_{2n}$ such that \[f_{i}(U w_{1}\ldots w_{i-1}w_{i}w_{i+1}\ldots w_{n-1}D)=U w_{1}\ldots w_{i-1}f(w_{i})w_{i+1}\ldots w_{n-1}D.\] We will denote  a unitary shift  by a vector of maps from $\mathfrak{D}_{2n}$ to itself of the form $(1_{1}, \ldots ,1_{i-1}, f_{i}, 1_{i+1}, \ldots, 1_{n-1})$, where $1_{k}$ is the identity map associated to the $i$-th coordinate.\\

An \textit{elementary shift} is a composition of unitary shifts. A \textit{shift path} of length $m$ $UWD \longrightarrow UW_1D\longrightarrow \cdots \longrightarrow UW_{m}D \longrightarrow UVD$ from $UWD$ to $UVD$ is a composition of $m$ elementary shifts. The set of all Dyck paths in a shift path between $UWD$ and $UVD$ will be denoted by $J(W,V)$. For notation, we introduce  the \textit{identity shift} as the elementary shift $(1_{1},\dots , 1_{n-1})$.\\

Suppose that a map $R: \mathfrak{D}_{2n}\rightarrow \mathfrak{D}_{2n}$ is defined by the application of successive elementary shifts to a given Dyck path. Then $R$ is said to be an  \textit{irreversible relation} over $\mathfrak{D}_{2n}$   if and only if  elementary shifts transforming Dyck paths (from one to the other) are not reversible. In other words, if an elementary shift $F=f_{p_{1}}\circ \dots \circ f_{p_{q}}$ transforms a Dyck path $UWD$ into a Dyck path $UVD$ then there is not an elementary shift $F'=f_{p_{1}}\circ \dots \circ f_{p_{q}}$ transforming $UVD$ into  $UWD$, for some $p,q \in \mathbb{Z}^{+}$.\\

If there exist two paths $G \circ F$ and $G' \circ F'$ of irreversible relations (of length 2) transforming a Dyck path $UWD$ into the Dyck path $UVD$ over $R$ in the following form:

	
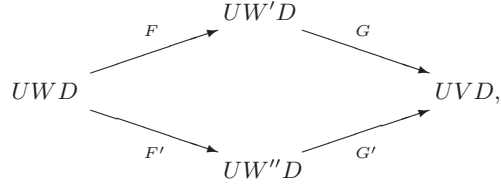
\begin{center}
	\begin{picture}(170,70)
		\put(0,30){$UWD$}
		\put(80,60){$UW'D$}
		\put(80,0){$UW''D$}
		\put(160,30){$UVD$,}
		\put(30,40){\vector(3,1){48}}
		\put(30,26){\vector(3,-1){48}}
		\put(110,10){\vector(3,1){48}}
		\put(110,56){\vector(3,-1){48}}
		\put(50,55){$_{F}$}
		\put(50,9){$_{F'}$}
		\put(130,55){$_{G}$}
		\put(130,9){$_{G'}$}
	\end{picture}
\end{center}

 with $W^{'}\neq W^{''}$. Then $G \circ F$  is said to be related with $G' \circ F'$ (denoted $G \circ F \sim_{R} G' \circ F'$) whenever $G'=F$ and $G=F'$.\par\bigskip

As for the case of diagonals \cite{Caldero}, Ca\~nadas and Rios \cite{Canadas} defined a  $\mathbb{F}$-linear  additive category $(\mathfrak{D}_{2n},R)$ based on Dyck paths, in this case,  \textit{objects} are $\mathbb{F}$-linear combinations of Dyck paths in $\mathfrak{D}_{2n}$ with  \textit{space of morphisms} from a  Dyck path $UWD$ to a Dyck path $UVD$ over $\mathbb{F}$ associated to $R$ being the vector space
\par\bigskip

\begin{centering}
 $\mathrm{Hom}_{(\mathfrak{D}_{2n},R)} (UWD,UVD)= \langle  \lbrace g \text{ } \vert \text{ }  g\hspace{0.1cm} \text{is a shift path associated to } R \rbrace \rangle / \langle  \sim_{R} \rangle$.\par\bigskip
 \end{centering}
  $\mathrm{Hom}_{(\mathfrak{D}_{2n},R)} (UWD,UVD)\neq0$ if and only if  there are shift paths transforming $UWD$ into $UVD$ and 

\begin{equation} \label{equation3.1}
\bigcap_{UW_{i}D \in J(W,V) } \text{Supp }UW_{i}D \neq \varnothing,\notag
\end{equation} 

for each shift path, with $UWD$ and $UVD$  in $\mathfrak{D}_{2n}$.\\

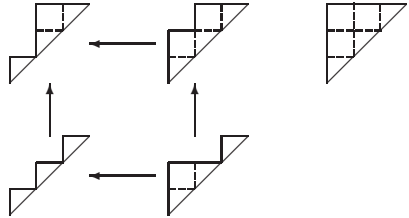
\begin{figure}[h!]
	\begin{center}
		\begin{picture}(170,80)
			\put(0,0){\line(1,1){30}}
			\put(0,0){\line(0,1){10}}
			\put(0,10){\line(1,0){10}}
			\put(10,10){\line(0,1){10}}
			\put(10,20){\line(1,0){10}}
			\put(20,20){\line(0,1){10}}
			\put(20,30){\line(1,0){10}}
			
			\put(0,50){\line(1,1){30}}
			\put(0,50){\line(0,1){10}}
			\put(0,60){\line(1,0){10}}
			\put(10,60){\line(0,1){20}}
			\put(10,80){\line(1,0){20}}
			\multiput(10,70)(3.6,0){3} {\line(1,0){2.6}} 
			\multiput(20,70)(0,3.6){3} {\line(0,1){2.6}} 
			
			\put(60,0){\line(1,1){30}}
			\put(60,0){\line(0,1){20}}
			\put(60,20){\line(1,0){20}}
			\put(80,20){\line(0,1){10}}
			\put(80,30){\line(1,0){10}}
			\multiput(60,10)(3.6,0){3} {\line(1,0){2.6}} 
			\multiput(70,10)(0,3.6){3} {\line(0,1){2.6}} 
			
			\put(60,50){\line(1,1){30}}
			\put(60,50){\line(0,1){20}}
			\put(60,70){\line(1,0){10}}
			\put(70,70){\line(0,1){10}}
			\put(70,80){\line(1,0){20}}
			\multiput(60,60)(3.6,0){3} {\line(1,0){2.6}} 
			\multiput(70,70)(3.6,0){3} {\line(1,0){2.6}} 
			\multiput(70,60)(0,3.6){3} {\line(0,1){2.6}} 
			\multiput(80,70)(0,3.6){3} {\line(0,1){2.6}} 
			
			\put(120,50){\line(1,1){30}}
			\put(120,50){\line(0,1){30}}
			\put(120,80){\line(1,0){30}}
			\multiput(120,60)(3.6,0){3} {\line(1,0){2.6}} 
			\multiput(120,70)(3.6,0){6} {\line(1,0){2.6}} 
			\multiput(130,60)(0,3.6){6} {\line(0,1){2.6}} 
			\multiput(140,70)(0,3.6){3} {\line(0,1){2.6}} 
			
			\put(15,30){\vector(0,1){20}}
			\put(70,30){\vector(0,1){20}}
			\put(55,15){\vector(-1,0){25}}
			\put(55,65){\vector(-1,0){25}}

		\end{picture}
	\end{center}
	
	\caption{ Elementary shifts in $(\mathfrak{D}_{6},R)$.}\label{figure3.1}
	
\end{figure}

Figure \ref{figure3.1} shows the elementary shifts over $(\mathfrak{D_{6}},R)$ associated to an irreversible relation $R$ defined over the set of all Dyck paths of length $6$. And such that,

\begin{equation}
R(UWD)=\begin{cases}
  f_{1}(UWD), & \mbox{ if }  w_{1}=UD,\\
  f_{2}(UWD),& \mbox{ if } w_{2}=UD.\notag
\end{cases}
\end{equation}

\begin{figure}[h!]
	\begin{center}
		
		\begin{picture}(160,15)
			
			\put(3,0){$Q=$}
			
			\put(30,3){\circle{5}} \put(60,3){\circle{5}} \put(90,3){\circle{5}} \put(120,3){\circle{5}} \put(150,3){\circle{5}}
			\put(62,3) {\vector(1,0){25}}
			\put(57.7,3){\vector(-1,0){25}}
			\put(122,3) {\vector(1,0){25}}
			\put(117.7,3){\vector(-1,0){25}}
			
			\put(28,8){$_{1}$}
			\put(58,8){$_{2}$}
			\put(88,8){$_{3}$}
			\put(118,8){$_{4}$}
			\put(148,8){$_{5}$}
			
		\end{picture}
		
	\end{center}
	
	\bigskip
	\bigskip
	
	\begin{center}
		
		\begin{picture}(276,228)
			
			\multiput(0,0)(96,0){3} {\line(1,1){36}}
			\multiput(48,48)(96,0){3} {\line(1,1){36}}
			\multiput(0,96)(96,0){3} {\line(1,1){36}}
			\multiput(48,144)(96,0){3} {\line(1,1){36}}
			\multiput(0,192)(96,0){3} {\line(1,1){36}}
			
			\multiput(0,0)(6,6){4} {\line(0,1){6}}
			\multiput(0,6)(6,6){4} {\line(1,0){6}} 
			\put(24,24){\line(0,1){12}} 
			\put(24,36){\line(1,0){12}} 
			
			\multiput(96,0)(6,6){2} {\line(0,1){6}}
			\multiput(96,6)(6,6){2} {\line(1,0){6}}
			\put(108,12){\line(0,1){12}} 
			\put(108,24){\line(1,0){6}} 
			\put(114,24){\line(0,1){6}}
			\put(114,30){\line(1,0){12}} 
			\put(126,30){\line(0,1){6}}
			\put(126,36){\line(1,0){6}}
			
			\put(192,0){\line(0,1){12}}
			\put(192,12){\line(1,0){6}}
			\put(198,12){\line(0,1){6}}
			\put(198,18){\line(1,0){12}}
			\multiput(210,18)(6,6){3} {\line(0,1){6}}  
			\multiput(210,24)(6,6){3} {\line(1,0){6}}
			
			\multiput(48,48)(6,6){2} {\line(0,1){6}} 
			\multiput(48,54)(6,6){2} {\line(1,0){6}}
			\put(60,60){\line(0,1){12}} 
			\multiput(60,72)(6,6){2} {\line(1,0){6}}
			\multiput(66,72)(6,6){2} {\line(0,1){6}}
			\put(72,84){\line(1,0){12}} 
			
			\put(144,48){\line(0,1){12}} 
			\multiput(144,60)(6,6){3} {\line(1,0){6}}
			\multiput(150,60)(6,6){3} {\line(0,1){6}} 
			\put(162,78){\line(1,0){12}}
			\put(174,78){\line(0,1){6}}
			\put(174,84){\line(1,0){6}}
			
			\put(240,48){\line(0,1){6}}
			\put(240,54){\line(1,0){6}}
			\put(246,54){\line(0,1){12}}
			\put(246,66){\line(1,0){12}} 
			\multiput(258,66)(6,6){3} {\line(0,1){6}}
			\multiput(258,72)(6,6){3} {\line(1,0){6}}
			
			\multiput(0,96)(6,6){2} {\line(0,1){6}}
			\multiput(0,102)(6,6){2} {\line(1,0){6}}
			\put(12,108){\line(0,1){12}}
			\put(12,120){\line(1,0){12}}
			\multiput(24,120)(6,6){2} {\line(0,1){6}}
			\multiput(24,126)(6,6){2} {\line(1,0){6}}
			
			\put(96,96){\line(0,1){12}}
			\multiput(96,108)(6,6){4} {\line(1,0){6}}
			\multiput(102,108)(6,6){4} {\line(0,1){6}}
			\put(120,132){\line(1,0){12}}
			
			\multiput(192,96)(30,30){2} {\line(0,1){6}}
			\multiput(192,102)(30,30){2} {\line(1,0){6}}
			\put(198,102){\line(0,1){12}}
			\multiput(198,114)(6,6){2} {\line(1,0){6}}
			\multiput(204,114)(6,6){2} {\line(0,1){6}}
			\put(210,126){\line(1,0){12}}
			
			\put(48,144){\line(0,1){12}}
			\multiput(48,156)(6,6){2} {\line(1,0){6}}
			\multiput(54,156)(6,6){2} {\line(0,1){6}}
			\put(60,168){\line(1,0){12}}
			\multiput(72,168)(6,6){2} {\line(0,1){6}}
			\multiput(72,174)(6,6){2} {\line(1,0){6}}
			
			\put(144,144){\line(0,1){6}}
			\put(144,150){\line(1,0){6}}
			\put(150,150){\line(0,1){12}}
			\multiput(150,162)(6,6){3} {\line(1,0){6}}
			\multiput(156,162)(6,6){3} {\line(0,1){6}}
			\put(168,180){\line(1,0){12}}
			
			\multiput(240,144)(6,6){3} {\line(0,1){6}}
			\multiput(240,150)(6,6){3} {\line(1,0){6}}
			\put(258,162){\line(0,1){12}}
			\put(258,174){\line(1,0){12}}
			\put(270,174){\line(0,1){6}}
			\put(270,180){\line(1,0){6}}
			
			\put(0,192){\line(0,1){12}}
			\put(0,204){\line(1,0){12}}
			\multiput(12,204)(6,6){4} {\line(0,1){6}}
			\multiput(12,210)(6,6){4} {\line(1,0){6}}
			
			\put(96,192){\line(0,1){6}}
			\put(96,198){\line(1,0){6}}
			\put(102,198){\line(0,1){12}}
			\put(102,210){\line(1,0){6}}
			\put(108,210){\line(0,1){6}}
			\put(108,216){\line(1,0){12}}
			\multiput(120,216)(6,6){2} {\line(0,1){6}}
			\multiput(120,222)(6,6){2} {\line(1,0){6}}
			
			\multiput(192,192)(6,6){3} {\line(0,1){6}}
			\multiput(192,198)(6,6){3} {\line(1,0){6}}
			\put(210,210){\line(0,1){12}}
			\put(210,222){\line(1,0){6}}
			\put(216,222){\line(0,1){6}}
			\put(216,228){\line(1,0){12}}
			
			\multiput(37,37)(96,0){3} {\vector(1,1){10}}
			\multiput(85,85)(96,0){2} {\vector(1,1){10}}
			\multiput(37,133)(96,0){3} {\vector(1,1){10}}
			\multiput(85,181)(96,0){2} {\vector(1,1){10}}
			\multiput(72,60)(96,0){2} {\vector(1,-1){30}}
			\multiput(24,108)(96,0){3} {\vector(1,-1){30}}
			\multiput(72,156)(96,0){2} {\vector(1,-1){30}}
			\multiput(24,204)(96,0){3} {\vector(1,-1){30}}

		\end{picture}

	\end{center}
	
	\caption{ Quiver $Q$ and the Dyck path version of the Auslander-Reiten quiver of $\mathrm{rep}\hspace{0.1cm}Q$ (see Corollary \ref{corol4.1}).}\label{figure3.3}
	
\end{figure}

If $\textbf{n}=\lbrace 1,2,\ldots, n\rbrace$ is an $n$-point chain then $\mathscr{C}_{\textbf{(1,n)}}$ stands for all \textit{admissible subchains} $\mathscr{C}$ of  $\textbf{n}$  with $\text{min }\mathscr{C}=1$ and $\text{max }\mathscr{C}=n$. Given an admissible subchain $\mathscr{C}=\{j_{1},\dots, j_{m},i_{1},\dots, i_{l}\}$ then two Dyck paths $D$ and $D'$ of length $2n$ are said to be related 
by a relation of type $R_{j_{1}\ldots j_{m}}^{i_{1} \ldots i_{l}}$ if there is an elementary shift associated to the points of the subchain transforming one into the other \cite{Canadas}.\\

We let $\mathfrak{C}_{2n}$ denote the subcategory of $(\mathfrak{D}_{2n},R_{j_{1}\ldots j_{m}}^{i_{1} \ldots i_{l}})$ whose objects are $\mathbb{F}$-linear combinations of Dyck paths with exactly $n-1$ peaks and related by a relation of type $R_{j_{1}\ldots j_{m}}^{i_{1} \ldots i_{l}}$.\\

The following results describe the structure of the category of Dyck paths \cite{Canadas}. We recall that a point $x$ where a Dyck path $UWD$ changes from the north to the east is said to be a \textit{peak} of the path. \\

 In the following theorem the symbol $Q$ is used to denote a quiver of type  $\mathbb{A}_{n-1}$ ($n\geq2$)  for which  $\lbrace i_{1},\dots, i_{l}\rbrace$ and $\lbrace j_{1},\dots, j_{m} \rbrace$ are the corresponding sets of sinks and sources. $\mathrm{rep}\hspace{0.1cm}Q$ denotes the corresponding category of representations.

\begin{teor}[\cite{Canadas}, Theorem 15]\label{teor3.1}

\textit{There is a categorical equivalence between categories $\mathfrak{C}_{2n}$ and $\mathrm{rep}\hspace{0.1cm} Q$.}

\end{teor} 
\addtocounter{corol}{2}
\begin{corol}[\cite{Canadas}, Corollary 16]\label{corol4.1}

\textit{There exists a bijection between the set of representatives of indecomposable representations  $\mathrm {Ind}\hspace{0.1cm}Q$ of $\mathrm{rep}\hspace{0.1cm}Q$ and the set of Dyck paths of length $2n$ with exactly $n-1$ peaks.}

\end{corol}

\begin{corol}[\cite{Canadas}, Corollary 17]\label{corol4.2}

\textit{The category $\mathfrak{C}_{2n}$ is an abelian category.}

\end{corol}

\section{Relationships Between Friezes and Dyck paths} \label{section4.1}

In this section, we present the main results of the paper. In particular, it is introduced the notion of integral diamond of Dynkin type $\mathbb{A}_{n}$, which are integer arrays used to build frieze patterns associated to triangulations of an $(n+3)$-polygon. 

\subsection{$\mathbb{A}_{n}$-Diamonds}

Let $\mathbf{R}$ be an integral domain then a \textit{diamond of Dynkin type} $\mathbb{A}_{n}$ or $\mathbb{A}_{n}$-diamond is an array $A=(a_{i, j})$ with entries in $\textbf{R}$, such entries satisfy conditions (D1) and (D2) associated to arrays with the following shape:

\[\xymatrix@=1mm{ & _{a_{2,0}} &  & & &\\
_{a_{1,1}} & & _{a_{2,1}}& & &\\
& _{a_{1,2}} & & _{\ddots} &\\ 
& & _{\ddots}& & _{a_{2,n-1}} &\\
&&&_{a_{1,n}} & & _{a_{2,n}}\\
&&&&_{a_{1,n+1}}&\\
}\]

\begin{enumerate}
\item [(D1)]$a_{2,0}=a_{1,n+1}=1$,
\item [(D2)]$a_{1,j}a_{2,j}-a_{2,j-1}a_{1,j+1}=1$ for $1\leq j \leq n$,
\end{enumerate}

where, $1$ is the identity element of $\mathbf{R}$. \par\bigskip

If $\mathbf{R}= \mathbb{Z}$ then $A$ is said to be a \textit{positive integral diamond of Dynkin type} $\mathbb{A}_{n}$, if it also satisfies the following condition (D3)
\begin{enumerate}
\item [(D3)]$a_{1,1}=a$ (or $a_{1,1}=a+m_{a}$), $a_{2,1}=a+m_{a}$ (or $a_{2,1}=a$) and $a_{1,2}=a^2+am_{a}-1$, with  $1 \leq a \leq \lfloor \frac{n+2}{2}\rfloor$,  $1 \leq m_{1}\leq n $ and $0 \leq m_{a} \leq n+2(1-a)$ if $a>1$.
\end{enumerate}
\par\bigskip
Henceforth, if it is not mentioned explicitly, $\mathbb{A}_{n}$-diamonds are assumed to be positive integral diamonds of Dynkin type $\mathbb{A}_{n}$.\\

Two  $\mathbb{A}_{n}$-diamonds $A$ and $B$ constitute a \textit{coupling}, denoted $A\models B$  if and only if $a_{2,j}=b_{1,j}$ for  $1\leq j\leq n$. \\

A set $\lbrace A^{t}\rbrace _{t\geq 0}$ is an $\mathbb{A}_{n}$-\textit{sequence of couplings} of $\mathbb{A}_{n}$ if and only if $A^{r}\models A^{r+1}$ for $r\geq 0$ ($X^0=X$ for any $\mathbb{A}_{n}$-diamond $X$). An $\mathbb{A}_{n}$-sequence $\lbrace A^{t}\rbrace _{t\geq 0}$ of couplings is a  $p$-\textit{cycle} if there exists $p \in \mathbb{N}$ such that $A^{t}=A^{t+p}$. If the $\mathbb{A}_{n}$-sequence of couplings $S_{t}=\lbrace A^{t}\rbrace _{t\geq 0}$ constitute a frieze pattern $\mathcal{F}$ then we will say that $S_{t}$ \textit{generates} $\mathcal{F}$. In this case, we point out that  Lemma \ref{lema4.1} proves that $A^0=A$ generates $S_{t}$, meaning that $S_{t}$ can be obtained from a sequence of couplings starting with $A$. \\

For example, let $\mathbf{R}=\mathbb{Z}$, the sets   $\lbrace A^{t}\rbrace _{t\geq 0}$ and $\lbrace B^{t}\rbrace _{t\geq 0}$ ($A^0$ and $B^0$ as shown below)  are  $\mathbb{A}_{1}$-sequences which are $2$-cycles with $A^{2k}=B^{2k+1}=A$, $A^{2k+1}=B^{2k}=B$ and  $k\geq 0$.

\begin{equation} \label{example1}
 \xymatrix@=1mm{ &&1& &&&&&1&\\
A=&1&&2 &&&B=&2&&1\\
&&1& &&&&&1&\\ 
}
\end{equation}

In general, it can be written an $\mathbb{A}_{n}$-sequence  $\lbrace A^{t}\rbrace _{t\geq 0}$ as an  $\mathbb{A}_{n}$-array  $C_{A^{t}}=(c_{i,j})$ such that $c_{t+1,j}=a^{t}_{1,j}$ and $c_{t+1,0}=c_{t+1,n+1}=1$, for $t\geq 0$. For the previous example, 

\[\xymatrix@=1mm{ &&1&&1&&1&\dots&&    &&&1&&1&&1&\dots&  \\
C_{A^{t}}=&1&&2&&1&&2&\dots&    &C_{B^{t}}= &2&&1&&2&&1&\dots \\
&&1&&1&&1&\dots&&  &&&1&&1&&1&\dots&  \\ 
}\]

$C_{A^{t}}$ and $ C_{B^{t}}$ are $\mathbb{A}_{1}$-arrays associated to $\lbrace A^{t}\rbrace _{t\geq 0}$ and $\lbrace B^{t}\rbrace _{t\geq 0}$, respectively.\\

If  the $\mathbb{A}_{n}$-sequence of couplings is  finite of length $m$, then  it can be associated to an infinity $\mathbb{A}_{n}$-array, $C^{m}_{A^{t}}=(c^{m}_{i,j})$ such that 
\begin{equation}\label{equation4.1.1}
 c^{m}_{(t+1)+km, j}=a^{t}_{1,j},\text{ } c^{m}_{(t+1)+km,0}=c^{m}_{(t+1)+km, n+1}=1,
\end{equation}

for $k \in \mathbb{Z}$. For any $\mathbb{A}_{n}$-sequence of couplings $\lbrace A^{t}\rbrace _{t\geq 0}$, it is possible to choose an $\mathbb{A}_{n}$-\textit{subsequence} $\lbrace B ^{z}\rbrace_{z \geq 0}$ with $B^{z}=A^{x+z}$ for some $x\geq t$. In particular, if $\lbrace A ^{t}\rbrace_{t\geq 0}$ is a $p$-cycle then the subsequence $\lbrace B ^{s_{0}}\rbrace_{0\leq s_{0} \leq p-1}$ such that $B^{s_{0}}=A^{t}$ is called the \textit{minimal} $p$-\textit{cycle} of $\lbrace A ^{t}\rbrace_{t\geq 0}$.\\

The following results give the main properties of diamonds of  Dynkin type $\mathbb{A}_{n}$.
\addtocounter{prop}{4}

\begin{prop}\label{prop4.1}
\textit{Let  $\lbrace A^{t}\rbrace_{t\geq 0}$ be a $p$-cycle and let  $B=\lbrace B ^{s_{0}}\rbrace_{0\leq s_{0} \leq p-1}$ be its minimal $p$-cycle. Then, the array $C^{p}_{B}$ is a frieze pattern of order $n+3$. In particular, $p$ divides $n+3$.}

\end{prop}

\textbf{Proof.} Let $C^{p}_{B}=(c^{p}_{ij})$ be the  infinity $\mathbb{A}_{n}$-array  associated to $\lbrace B ^{s_{0}}\rbrace_{0\leq s_{0} \leq p-1}$,  identity (\ref{equation4.1.1}) implies that
\[ c^{p}_{(s_{0}+1)+kp, j}=a^{s_{0}}_{1,j}, \quad  c^{p}_{(s_{0}+1)+kp, 0}=c^{p}_{(s_{0}+1)+kp, n+1}=1,\]
for $k \in \mathbb{Z}$, since that $\lbrace A^{t}\rbrace_{t\geq 0}$ is a $p$-cycle, then $C^{p}_{B}$ is a frieze pattern. $\hfill\square$

\begin{prop}\label{prop4.2}
\textit{Let $\lbrace A^{t}\rbrace_{t\geq 0}$ be a $p$-cycle of length $2p$, then the subsequences $\lbrace B^{s_{i}}\rbrace$ with ${0\leq s_{i} \leq p-1}$  generate the same frieze pattern of order $n+3$, for $0\leq i \leq p-1$, and $B^{s_{i}}=A^{i+s_{i}}$.}

\end{prop}

\textbf{Proof.} Let $\lbrace A^{t}\rbrace_{t\geq 0}$ be a $p$-cycle of length $2p$, let $C^{p}_{A}=(c^{p}_{ij})$ and $C^{p}_{B}=(c^{p'}_{ij})$ be the infinity arrays of the subsequences $A=\lbrace B^{s_{i}}\rbrace_{0 \leq s_{i} \leq p-1}$ and $B=\lbrace B_{s_{i'}}\rbrace_{0 \leq s_{i'} \leq p-1}$  of $\lbrace A^{t}\rbrace_{t\geq 0}$  for $0 \leq i<i'\leq p-1$. The following identities (\ref{idp4}) hold by applying  the translation $s_{i'}=s_{i}-|i'-i|$, 
\begin{equation}\label{idp4}
\begin{split}
c^{p}_{s_{i'}+1+kp, j}=a_{1j}^{s_{i'}+i'}= a_{1j}^{s_{i}-|i'-i|+i'}=a_{ij}^{s_{i}+i}=c^{p}_{s_{i}+1+kp, j}. 
\end{split}
\end{equation}
We are done.  $\hfill\square$\\

\addtocounter{lema}{6}
\begin{lema}\label{lema4.1}
\textit{Let  $\lbrace A^{t}\rbrace_{t\geq 0}$ be a sequence of couplings. Then  $\lbrace A^{t}\rbrace_{t\geq 0}$ is generated by  $A^{0}$. In particular,  $A^{0}$ generates a $p$-cycle for some $p>0$.}
\end{lema}

\textbf{Proof.} Let $\lbrace A^{t}\rbrace_{t\geq 0}$ be a sequence, then 
\begin{equation}
a^{x}_{2,j}=\frac{1+(a^{x}_{2,j-1})(a^{x-1}_{2,j+1})}{a^{x-1}_{2,j}},\notag
\end{equation}

for $1 \leq j \leq n$, and $x \geq t$, $a^{x}_{2,j}$ can be written  by using the set $\{a^{0}_{2,j} \}_{1\leq j \leq n}$ for $x>0$. In particular, the set $\lbrace a^{0}_{2,1}, \dots, a^{0}_{2,n}\rbrace$ is a seed  of  the cluster algebra  associated to  the linearly oriented quiver of type $\mathbb{A}_{n}$.  Since the cluster variables  are finite in the case $\mathbb{A}_{n}$, then there is $p=n+3$ (in some cases, it is not minimal) such that $A^0=A^{p}$. $\hfill\square$

\addtocounter{teor}{5}
\begin{teor}\label{teor4.1}
\textit{Let $A$ be a diamond of Dynkin type $\mathbb{A}_{n}$ then $A$  generates a frieze pattern.}

\end{teor}

\textbf{Proof.} It is a direct consequence of Lemma \ref{lema4.1}, and Proposition \ref{prop4.1}.$\hfill\square$ \\

For instance, diamonds $A$ and $B$ given in (\ref{example1}) generate the following  frieze pattern.

\[\xymatrix@=1mm{ \ldots &1&&1& & 1 & & 1 && 1&\ldots\\
 &\ldots &1& &2& &1& &2&\ldots&\\
\ldots &1&&1& &1& &1& &1&\ldots\\ 
}\]

\subsection{Seed Vectors}

In this section, we give an algorithm to build a family of positive integral frieze vectors associated to  the linearly oriented quiver of type $\mathbb{A}_{n}$.  These vectors  allow  to find out a connection between the positive integral diamonds  of Dynkin type $\mathbb{A}_{n}$, triangulations, and Dyck paths. \\

Let $A$ be a diamond of Dynkin type $\mathbb{A}_{n}$ then we  can write its first column as a vector with the form $v_{A}=(a_{1}, \dots, a_{n})$, where $a_{j}=a_{1,j}$. In such a case, we say that $v_{A}$ is associated to $A$ and that $v_{A}$ generates $A$.

\addtocounter{prop}{2}
\begin{prop}\label{prop4.3}
\textit{If $v=(a_{1},\dots, a_{n})$ is a vector associated to a positive integral diamond  of Dynkin type $\mathbb{A}_{n}$ with $a_{n}=1$, then the vector $v'=(a_{1}, \dots, a_{i},a_{i}+a_{i+1},a_{i+1},\dots ,a_{n-1})$ is also associated to a positive integral diamond of Dynkin type $\mathbb{A}_{n}$, for $1 \leq i<n$.}
\end{prop}

\textbf{Proof.} Let $v_{A}=(a_{1},\dots, a_{n})$ be a vector associated to a positive integral diamond $A=(a_{j, m})$ of Dynkin type $\mathbb{A}_{n}$, then we take  the  vector $v_{A+i}=(a_{1}, \dots, a_{i},a_{i}+a_{i+1},a_{i+1},\dots ,a_{n-1})$ and  the array $A+i$ of the following form:  

\begin{equation}
b_{1,m}=\begin{cases}
  a_{1,m}, & \mbox{if $m \leq i$,}\\
 a_{1,i}+a_{1,i+1}, & \mbox{if $m= i+1$,}\\\notag
 a_{1,m-1}, & \mbox{if $m> i+1$,}
\end{cases}
\end{equation}

and 
\begin{equation}
b_{2,m}=\begin{cases}
  a_{2,m}, & \mbox{if $m \leq i-1$,}\\
  a_{2,i-1}+a_{2,i}, & \mbox{if $m= i$,}\\\notag
a_{2,m-1}, & \mbox{if $m \geq i+1$,}
\end{cases}
\end{equation}
 
 then $b_{1,m}b_{2,m}-b_{2,m-1}b_{2,m+1}=1$, for $1 \leq m \leq n$ and $1 \leq i <n$. Therefore $A+i$ is a positive integral diamond of Dynkin type $\mathbb{A}_{n}$. $\hfill\square$

\begin{prop}\label{prop4.4}
\textit{For each vector $v_{n,z}=(a_{1},\dots ,a_{n})$ with}

\begin{equation} \label{equation4.2.8}
 a_{i}=\begin{cases}
  z+1-i, & \mbox{if $i<z$,}\\
 1, & \mbox{if $i \geq z$.}
\end{cases}
\end{equation}

 \textit{there is associated a unique positive integral diamond of Dynkin type $\mathbb{A}_{n}$, for $z \in \{1,\dots,n+1 \}$.}
 
\end{prop}

\textbf{Proof.} Let $v_{n,z}$ be a vector and let $z$ be a natural number between $1$ and $n+1$, we define a positive integral diamond $A$  with $a_{1,i}=a_{i}$ and $a_{2,i}=b_{i}$  where

\begin{equation} \label{equation4.1.2}
b_{i}=\begin{cases}
  1, & \mbox{if $i<z$,}\\
 i+2-z, & \mbox{if $i \geq z$,}
\end{cases}
\end{equation}

then $a_{1,i}a_{2,i}-a_{2,i-1}a_{2,i+1}=1$ for $1 \leq i \leq n$. $\hfill\square$

\addtocounter{Nota}{10}
\begin{Nota}

$v_{n,z}$  is called a \textit{seed vector}.  The vector $v^{n,z}=(b_{1},\dots , b_{n})$  defines a positive integral diamond $B$ of  Dynkin type $\mathbb{A}_{n}$ such that $b_{2,i}$ satisfies the following identity

\begin{equation} \label{equation4.1.3}
b_{2,i}=\begin{cases}
  i-1, & \mbox{if $i<z-1$,}\\
 (b_{1,i}+1)z-1, & \mbox{if $z-2< i< n$,}\\\notag
 z, & \mbox{if $i=n$.}
\end{cases}
\end{equation}
and  $b_{i}=b_{1,i}$ is defined as in (\ref{equation4.1.2}).
\end{Nota}

\addtocounter{prop}{1}

\begin{prop}
\textit{The positive integral diamonds $A$ and $B$ of Dynkin type $\mathbb{A}_{n}$ generated by vectors $v_{n, z}$ and $v^{n, z}$ respectively constitute a coupling.} 
\end{prop}
\textbf{Proof.} It is a direct consequence of  Proposition \ref{prop4.4} and Lemma \ref{lema4.1}.$\hfill\square$\\

The number of ways of applying recursively Proposition \ref{prop4.3} to a  vector $w_{A}=(a_{1},\dots,$ $a_{z-1},1, \dots, 1) \in \mathbb{N}^{n}$  is given by the next identity (denoted by $f_{n,z}$),

\begin{equation}
f_{n, z}=\begin{cases}
  \sum _{i=z-1}^{n} f_{n-1,i}, & \mbox{if $z >1$,}\\\notag
 \sum _{i=1}^{n} f_{n-1,i}, & \mbox{if $z=1$,}
\end{cases}
\end{equation}

where it is included the trivial move $w_{A+0}=w_{A}$, for $n>1$, and any $z \in \lbrace1, \dots, n+1\rbrace$. Actually, these numbers can be arranged  as follows:
\begin{equation} \label{fibt}
\xymatrix @=0.1mm{f_{1,2}& f_{1,1}& & \\
    f_{2,3}&f_{2,2}&f_{2,1}&\\
    f_{3,4}&f_{3,3}&f_{3,2}&f_{3,1}\\
    f_{4,5}&f_{4,4}&f_{4,3}&f_{4,2}& f_{4,1}\\
    \vdots &&&& & \ddots }
\end{equation}  

for any of the vectors $w_{A}$. Since the first  choices are $v_{1,1}=(1)$ and $v_{1,2}=(2)$, then $f_{1,1}=1$ and $f_{1,2}=1$. The previous triangle (\ref{fibt}) appears in the On-Line Encyclopedia of Integer Sequences (OEIS) as A009766 (Catalan triangle \cite{OEIS9}). In particular, we generate all positive integral diamonds of Dynkin type $\mathbb{A}_{n}$ via the seed vectors $v_{n, z}$. For example,  for $n=3$, all vectors that generate positive integral diamonds of Dynkin type $\mathbb{A}_{3}$ are:
 
\[\xymatrix @=1mm{(1,1,1)&(1,1,2)&(1,2,1)&(1,2,3)&(1,3,2)& (2,1,1)&(2,1,2)\\ 
    (2,3,1)&(2,3,4)&(2,5,3)& (3,2,1)&(3,2,3)&(3,5,2)&(4,3,2)\\}\]

\begin{figure}[h!]
	\begin{center}
		\begin{picture}(94.5,94.5)
			\put(0,0){\line(0,1){52.5}}
			\put(0,52.5){\line(1,0){31.5}}
			\put(31.5,52.5){\line(0,1){10.5}}
			\put(31.5,63){\line(1,0){10.5}}
			\put(42,63){\line(0,1){31.5}}
			\put(42,94.5){\line(1,0){52.5}}
			\put(0,0){\line(1,1){94.5}}
			\multiput(0,10.5)(3.6,0){3} {\line(1,0){2.6}} 
			\multiput(0,21)(3.6,0){6} {\line(1,0){2.6}} 
			\multiput(0,31.5)(3.6,0){9} {\line(1,0){2.6}}
			\multiput(0,42)(3.6,0){12} {\line(1,0){2.6}}
			\multiput(31.5,52.5)(3.6,0){6} {\line(1,0){2.6}}
			\multiput(42,63)(3.6,0){6} {\line(1,0){2.6}}
			\multiput(42,73.5)(3.6,0){9} {\line(1,0){2.6}}
			\multiput(42,84)(3.6,0){12} {\line(1,0){2.6}}
			\multiput(10.5,10.5)(0,3.6){12} {\line(0,1){2.6}}
			\multiput(21,21)(0,3.6){9} {\line(0,1){2.6}}
			\multiput(31.5,31.5)(0,3.6){6} {\line(0,1){2.6}}
			\multiput(42,42)(0,3.6){6} {\line(0,1){2.6}}
			\multiput(52.5,52.5)(0,3.6){12} {\line(0,1){2.6}}
			\multiput(63,63)(0,3.6){9} {\line(0,1){2.6}}
			\multiput(73.5,73.5)(0,3.6){6} {\line(0,1){2.6}}
			\multiput(84,84)(0,3.6){3} {\line(0,1){2.6}}
		\end{picture}
	\end{center}
	\caption{A Dyck path of length  $18$ associated to the vector $v_{G}=(5,4,3,3,5,4,3,2)$.}\label{figure4.1}
\end{figure}
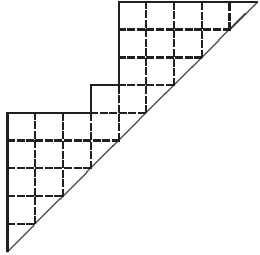

Let $G= UD \dots UD \dots$ be a Dyck path of length $2n$ and  let $m_{i}$ be the number of $U$s before the occurrence of the $i$-th $D$ in $G$ then $G$ can be written  as a vector $v_{G}=(v_{1},\dots, v_{n-1})$ where $v_{i}=m_{i}-i+1$. As an example, consider the Dyck path $G$ shown in Figure \ref{figure4.1}, which has associated the vector $v_{G}=(5,4,3,3,5,4,3,2)$.\\

In what follows, it is defined a map $T_{i}$ between vectors associated to positive integral diamonds of Dynkin type $\mathbb{A}_{n}$ and Dyck paths by using a relation over the coordinates of a vector $u=(a_{1},\ldots, a_{m})$. $T_{i}$ is defined in such a way that, $T_{i}:\mathbb{N}^m  \rightarrow  \mathbb{N}$ and: 

\begin{itemize}

\item If $a_{i}-a_{l}>0$ for some $l \in \lbrace 1, \ldots, i \rbrace$, then we select the maximum of such indexes $\mathrm{max}\hspace{0.1 cm}\lbrace l\rbrace $ and write $r_{1}=a_{i}-a_{l}$. In the same way, it is chosen $\mathrm{max}\hspace{0.1 cm}\lbrace l\rbrace $ such that $r_{1}-a_{l}> 0$ and write $r_{2}=r_{1}-a_{l}$, this process ends if $\{l\mid r_{t}-a_{l}> 0\}=\varnothing$. Thus, if $u\in\mathbb{N}^{m}$ then $T_{i}(u)=r_{t} + t$, for some $t$.

\item If $a_{i}-a_{l}\leq 0$ for all $l  \in \lbrace 1, \ldots, i \rbrace$, then $T_{i}(u)= a_{i}$.

\end{itemize}

For instance, for the vector $u=(14,52,4,23,9,2)$,  it holds that $T_{1}(u)=14$, $T_{2}(u)=13$, $T_{3}(u)=4$, $T_{4}(u)=8$, $T_{5}(u)=3$, and $T_{6}(u)=2$.

\begin{prop} \label{prop4.10}
\textit{Let $v_{n,z}$ be a seed vector, then $(T_{1}(v_{n,z}), \dots, T_{n}(v_{n,z}))$ defines a Dyck path of length $2(n+1)$.}
\end{prop}

\textbf{Proof.} For any $z \in \lbrace  1, \dots, n+1 \rbrace$, $T_{i}(v_{n,z})=a_{i}$ with $a_{i}$ given by identity (\ref{equation4.2.8}), there is a word  $G_{v_{n,z}}=w_{1}\dots w_{2(n+1)}$  in the free monoid $ \lbrace U,D\rbrace^{\ast}$ such that

\begin{equation}
G_{v_{n,z}}= \underbrace{U \dots U}_{z-1} \underbrace{D \dots D}_{z-1}UDUD \dots UDUD,\notag
\end{equation}
for any left factor $u_{s}$ in $G_{v_{n,z}}$  of length $s \in \lbrace 1 , \dots, 2(n+1)\rbrace$,  $0 \leq |u_{s}|_{U}- |u_{s}|_{D} \leq z-1$, therefore $G_{v_{n,z}} \in \mathfrak{D}_{2(n+1)}$.$\hfill\square$\\

\begin{prop} \label{prop4.11}
\textit{Let $v_{A}=(a_{1}, \dots, a_{n})$ be a vector associated to a positive integral diamond $A$ of Dynkin type $\mathbb{A}_{n}$ with $a_{n}=1$, such that $(T_{1}(v_{A}), \dots, T_{n}(v_{A}))$ defines a Dyck path in $\mathfrak{D}_{2(n+1)}$. Then $(T_{1}(v_{A+i}), \dots, T_{n}(v_{A+i}))$ also defines a Dyck path in $\mathfrak{D}_{2(n+1)}$.} 
\end{prop}

\textbf{Proof.} Let $v_{A}=(a_{1}, \dots, a_{n})$ be a vector associated to a positive integral diamond $A$ with $a_{n}=1$, then there exists a Dyck path $G_{v_{A}} \in \mathfrak{D}_{2(n+1)}$ such that any left factor $u_{s}$ of length $s$ satisfies $|u_{s}|_{U} \geq |u_{s}|_{D}$ for  $1 \leq s \leq 2(n+1)$. Let $v_{A+i}$ be a vector associated to the positive integral diamond $A+i$ with
\begin{equation}
T_{m}(v_{A+i})=\begin{cases}
  T_{m}(v_{A}), & \mbox{if $1 \leq m \leq i$,}\\
 T_{m}(v_{A})+1, & \mbox{if $m= i+1$,}\\\notag
 T_{m-1}(v_{A}), & \mbox{if $m\geq i+1$,}
\end{cases}
\end{equation}
then there is a word $G_{A+i}=w'_{1}, \dots, w'_{2(n+1)}$ in $\lbrace U,D\rbrace^{\ast}$, where we take the index $m_{1}$   of the $i$-th $D$  in $G_{A+i}$,  any left factor $u'_{s}$ in $G_{A+i}$ satisfies the identities

\begin{equation}
|u'_{s}|_{U}=\begin{cases}
  |u_{s}|_{U}, & \mbox{if $1 \leq s \leq m_{1}$,}\\
 |u_{m_{1}}|_{U}+1, & \mbox{if $s= m_{1}+1$,}\\\notag
 |u_{s-2}|_{U}+1, & \mbox{if $s \geq m_{1}+2$,}
\end{cases} 
\end{equation}
and 
\begin{equation}
|u'_{s}|_{D}=\begin{cases}
  |u_{s}|_{D}, & \mbox{if $1 \leq s \leq m_{1}$,}\\
 |u_{m_{1}}|_{D}, & \mbox{if $s= m_{1}+1$,}\\\notag
 |u_{s-2}|_{U}+1, & \mbox{if $s \geq m_{1}+2$,}
\end{cases}
\end{equation}
then, we have the following possibilities:

\begin{itemize}
\item If $1 \leq s \leq m_{1}$, $|u'_{s}|_{U}=|u_{s}|_{U} \geq |u_{s}|_{D}=|u'_{s}|_{D}$.
\item If $s=m_{1}+1$, $|u'_{m_{1}+1}|_{U}=|u_{m_{1}}|_{U}+1 > |u_{m_{1}}|_{D}=|u'_{m_{1}+1}|_{D}$.
\item If $m_{1}+2 \leq s \leq 2(n+1)$, $|u'_{s}|_{U}=|u_{s-2}|_{U}+1 \geq |u_{s-2}|_{D}+1=|u'_{s}|_{D}$.
\end{itemize}

Therefore, $G_{A+i} \in \mathfrak{D}_{2(n+1)}$.$\hfill\square$\\
\addtocounter{lema}{7}

\begin{lema}\label{lema4.2}
\textit{There is a bijective correspondence between  the set of all  vectors associated to positive integral diamonds of Dynkin type  $\mathbb{A}_{n}$ and the set of all Dyck paths  of length  $2(n+1)$.}
\end{lema}

\textbf{Proof.} Let $\mathbb{D}_{\mathbb{A}_{n}}$  be the set of all  vectors associated to positive integral diamonds of Dynkin type  $\mathbb{A}_{n}$ and  let $\mathfrak{D}_{2(n+1)}$ be  the set of all Dyck paths  of length  $2(n+1)$ then we define a map $f:\mathbb{D}_{\mathbb{A}_{n}}  \rightarrow \mathfrak{D}_{2(n+1)} $ with $f(u_{A})=\big (T_{1}(u_{A}), \ldots , T_{n}(u_{A})\big )$, Propositions \ref{prop4.10} and \ref{prop4.11} allow us to establish that $f$ is well defined. In order to prove that the map $f$ is injective, suppose that $u_{A}\neq v_{B}$, we take the minimum $l$ such that $u_{l}\neq v_{l}$. If $l=1$ then $T_{1}(u_{A})\neq T_{1}(v_{B})$. If $l>1$,   $u_{l}=m(u_{l-1})+ a$ and $v_{l}=m'(u_{l-1})+a$  with  $m\neq m'$ is a consequence of Proposition \ref{prop4.3} then $r_{t_{u_{l}}}\neq r_{t_{v_{l}}}$, therefore $T_{l}(u_{A}) \neq T_{k}(v_{B})$. $\hfill\square$\\

An alternative way of writing a Dyck path $G \in \mathfrak{D}_{2(n+1)}$ can be defined by using a vector $\lambda_{G}=(\lambda_{1},\dots, \lambda_{n})$ where $\lambda_{i}$ is the number of $D$s  before the occurrence of the $(n+2-i)$-th $U$ in $G$. For example, Dyck path in Figure \ref{figure4.1} has associated the following vector $\lambda_{G}=(4,4,4,3,0,0,0,0)$.\\
     
Let $\lambda$ be a vector associated to a Dyck path of length $2(n+1)$ then a triangulation of  an $(n+3)$-polygon can be realized by $\lambda$ as follows:
 
 \begin{itemize}
\item Fix a labeling  for the vertices of polygon  $K_{0}^{n+3}=(v^{n+3}_{0}, \dots, v^{n+3}_{n+2})$  with $v^{n+3}_{i}=i$, for $0 \leq i \leq n+2$. 
\item For $\lambda _{i}$, we draw a diagonal $l_{i}^{\lambda _{i}}$ between $\lambda_{i}$ and $\lambda_{i}+2$. Afterwards, we label the last polygon with $n+3-i$ vertices $K_{i}^{n+3-i}=(v^{n+3-i}_{0}, \dots, v^{n+3-i}_{n+2-i})$, and

\begin{equation*}
v^{n+3-i}_{j}=\begin{cases}
 v^{n+3-(i-1)}_{j}, & \mbox{if $j \leq \lambda_{i}$,}\\
 v^{n+3-(i-1)}_{j+1}-1, & \mbox{if $j> \lambda_{i}$,}
\end{cases}
\end{equation*}
\end{itemize}
 for $i=1, \dots ,n$.\\
 
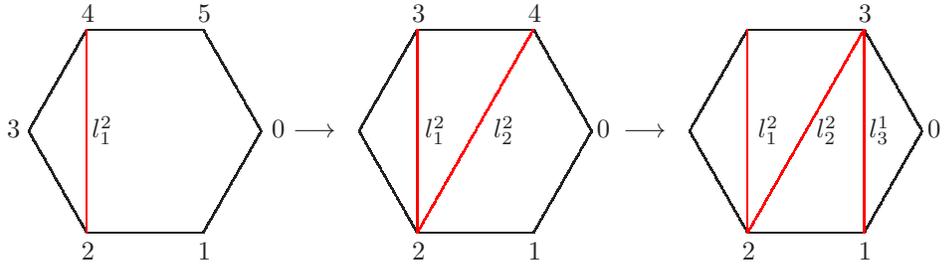
\begin{figure}[h!]
\begin{center}
	\setlength{\unitlength}{2pt}
\begin{picture}(160,42)

\qbezier(8,0)(19,0)(30,0)
\qbezier(30,0)(35.5,9.625)(41,19.25)
\qbezier(7,1.75)(7.5,0.875)(8,0)
\qbezier(-3,19.25)(2.5,9.625)(8,0)
\qbezier(-3,19.25)(2.5,28.875)(8,38.5)
\qbezier(8,38.5)(20,38.5)(30,38.5)
\qbezier(30,38.5)(35.5,28.875)(41,19.25)
\color{red}
\qbezier(8,0)(8,19.25)(8,38.5)

\normalcolor


\qbezier(70.5,0)(81.5,0)(92.5,0)
\qbezier(92.5,0)(98,9.625)(103.5,19.25)
\qbezier(69.5,1.75)(70,0.875)(70.5,0)
\qbezier(59.5,19.25)(65,9.625)(70.5,0)
\qbezier(59.5,19.25)(65,28.875)(70.5,38.5)
\qbezier(70.5,38.5)(82.5,38.5)(92.5,38.5)
\qbezier(92.5,38.5)(98,28.875)(103.5,19.25)
\color{red}
\qbezier(70.5,0)(70.5,19.25)(70.5,38.5)
\qbezier(70.5,0)(81.5,19.25)(92.5,38.5)

\normalcolor


\qbezier(133,0)(144,0)(155,0)
\qbezier(155,0)(160.5,9.625)(166,19.25)
\qbezier(132,1.75)(132.5,0.875)(133,0)
\qbezier(122,19.25)(127.5,9.625)(133,0)
\qbezier(122,19.25)(127.5,28.875)(133,38.5)
\qbezier(133,38.5)(145,38.5)(155,38.5)
\qbezier(155,38.5)(160.5,28.875)(166,19.25)
\color{red}
\qbezier(133,0)(133,19.25)(133,38.5)
\qbezier(133,0)(144,19.25)(155,38.5)
\qbezier(155,38.5)(155,19.25)(155,0)

\normalcolor

\put(7,-5){$2$}
\put(29,-5){$1$}
\put(29,40){$5$}
\put(7,40){$4$}
\put(43,18){$0$}
\put(-7,18){$3$}

\put(9,18){$l_{1}^2$}

\put(69.5,-5){$2$}
\put(91.5,-5){$1$}
\put(91.5,40){$4$}
\put(69.5,40){$3$}
\put(104.5,18){$0$}

\put(72,18){$l_{1}^2$}
\put(85,18){$l_{2}^2$}

\put(132,-5){$2$}
\put(154,-5){$1$}
\put(154,40){$3$}
\put(167,18){$0$}

\put(135,18){$l_{1}^2$}
\put(146,18){$l_{2}^2$}
\put(156,18){$l_{3}^1$}

\put(47,18){$\longrightarrow$}
\put(109.5,18){$\longrightarrow$}

\end{picture}
 \end{center}
 \caption{ Example of a triangulation realized by a Dyck path of length 8.}\label{figure4.3}
\end{figure}

The previous algorithm describes that  if $l_{i}^{\lambda_{i}}$ is a diagonal then it does not cross the diagonals $l_{1}^{\lambda_{1}},\dots,l_{i-1}^{\lambda_{i-1}}$ for $1 \leq i\leq n$. For instance, let $\lambda_{G}=(2,2,1)$ be the vector associated to $G=UDUDUUDD$, then the triangulation of $\lambda_{G}$ is shown in Figure \ref{figure4.3}.\\

If we fix a labeling $K$ over all vertices of a polygon with $n+3$ vertices, a triangulation $T$ is written as a sequence $T=(l_{1}^{v_{1}},\dots, l_{n}^{v_{n}})$, where $v_{i}$ belongs to the set of vertices.

\begin{lema}\label{lema4.3}

\textit{There is a bijective correspondence between the set of all triangulations of a polygon with $n+3$ vertices and the set of all Dyck paths  of length  $2(n+1)$.}

\end{lema}

\textbf{Proof.}  Let  $\mathcal{T}_{n}$ be the set of all triangulations of a polygon with $n+3$ vertices then we can define a map $g:\mathfrak{D}_{2(n+1)} \rightarrow \mathcal{T}_{n}$ with $g(\lambda)=T_{\lambda}$. In order to prove that $g$ is one to one, we fix a labeling $K$ and suppose that $g(\lambda_{G})=g(\sigma_{G'})$, then $(l_{1}^{\lambda_{1}},\dots, l_{n}^{\lambda_{n}})=(l_{1}^{\sigma_{1}}, \dots, l_{n}^{\sigma_{n}})$, provided that $l_{j}^{\lambda_{j}}=l_{j}^{\sigma_{j}}$. Since by definition there are diagonals connecting vertices $\lambda_{j}$ with $(\lambda_{j} +2)$ and $\sigma_{j}$ with $(\sigma_{j}+2)$, therefore $\lambda_{j}= \sigma_{j}$ for $j=1,\dots, n$. We are done.$\hfill\square$\\

The next theorem presents the main result regarding relationships between positive integral diamonds of Dynkin type $\mathbb{A}_{n}$ and triangulations of an $(n+3)$-polygon.
\addtocounter{teor}{8}

\begin{teor}\label{teor4.2}
\textit{There is a bijective correspondence between the set of all  vectors associated to positive integral diamonds of Dynkin type  $\mathbb{A}_{n}$  and triangulations  of a polygon with $n+3$ vertices.}
\end{teor}

\textbf{Proof.} We fix a labeling  $K$ in a polygon with $n+3$ vertices, the map $F:\mathbb{D}_{\mathbb{A}_{n}} \rightarrow \mathcal{T}_{n}$ defined by the formula
\begin{equation}\label{F}
\begin{split}
F(v_{A})&=(g \circ f)(v_{A})
\end{split}
\end{equation}

is a bijection (see Lemmas \ref{lema4.2} and \ref{lema4.3}). $\hfill\square$\\

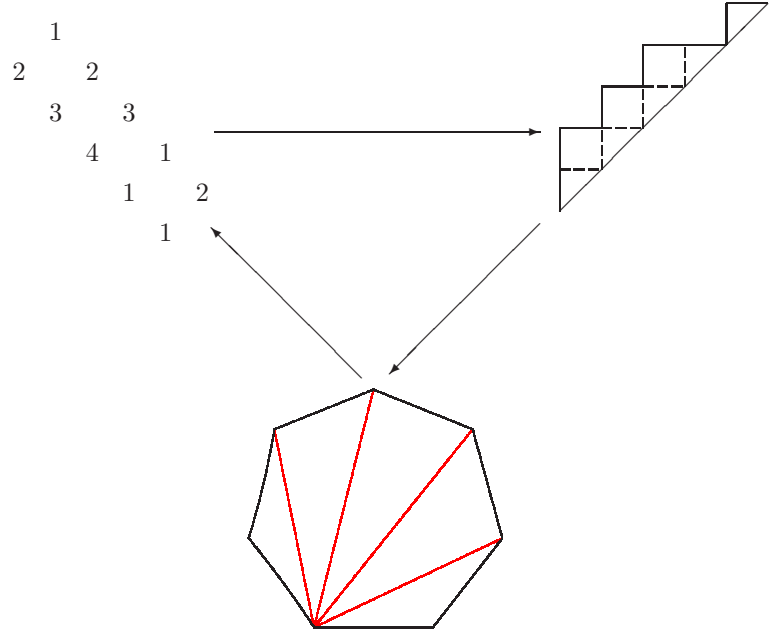
\begin{figure}[h!]
	\begin{center}
		\setlength{\unitlength}{1.5pt}
		\begin{picture}(100,150)
			
			\put(-45,150){{{{\xymatrix@=1mm{ &1&&&& \\
								2&&2&&&\\
								&3&&3&&\\
								&&4&&1&\\
								&&&1&&2\\
								&&&&1&\\ 
			}}}}}	
			\put(8,125){\vector(1,0){82}}
			\put(90,102){\vector(-1,-1){38}}
			\put(45,63){\vector(-1,1){38}}
			

			\put(95,105){\line(1,1){52.5}}
			\put(95,105){\line(0,1){21}}
			\put(95,126){\line(1,0){10.5}}
			\put(105.5,126){\line(0,1){10.5}}
			\put(105.5,136.5){\line(1,0){10.5}}
			\put(116,136.5){\line(0,1){10.5}}
			\put(116,147){\line(1,0){21}}
			
			\put(137,147){\line(0,1){10.5}}
			\put(137,157.5){\line(1,0){10.5}}

			\multiput(95,115.5)(3.6,0){3} {\line(1,0){2.6}} 
			\multiput(105.5,126)(3.6,0){3} {\line(1,0){2.6}} 
			\multiput(116,136.5)(3.6,0){3} {\line(1,0){2.6}}

			\multiput(105.5,115.5)(0,3.6){3} {\line(0,1){2.6}}
			\multiput(116,126)(0,3.6){3} {\line(0,1){2.6}}
			\multiput(126.5,136.5)(0,3.6){3} {\line(0,1){2.6}}
			
			\color{red}
			\qbezier(33,0)(56.75,11.25)(80.5,22.5)
			\qbezier(33,0)(53,25)(73,50)
			\qbezier(33,0)(40.5,30)(48,60)
			\qbezier(33,0)(28,25)(23,50)
			\normalcolor
			\qbezier(33,0)(48,0)(63,0)
			\qbezier(63,0)(71.75,11.25)(80.5,22.5)
			\qbezier(80.5,22.5)(76.75,36.25)(73,50)
			\qbezier(73,50)(60.5,55)(48,60)
			\qbezier(48,60)(35.5,55)(23,50)

			\qbezier(23,50)(20.75,36.25)(16.5,22.5)
			\qbezier(16.5,22.5)(25.75,11.25)(33,0)
			
		\end{picture}
	\end{center}
	\caption{Map $F$ defined in the proof of Theorem \ref{teor4.2} allows to establish identifications between an $\mathbb{A}_{4}$-diamond, a Dyck path of length 10 and a triangulation of a heptagon.}\label{figure4.4}
\end{figure}

Figure \ref{figure4.4} presents an example of the bijective correspondence between a positive integral diamond of Dynkin type $\mathbb{A}_{4}$, a Dyck path of length $10$, and a triangulation of a polygon with $7$ vertices.

\subsection{Frieze Patterns and Dyck Paths}

In this section, we describe an algebraic interpretation of frieze patterns as a direct sum of indecomposable objects of Dyck paths categories.

\addtocounter{lema}{1}
\begin{lema}\label{lema4.4}
\textit{Vectors $v_{n,z}$ and $v^{n,z}$ realize the same triangulation except for one anti-clockwise rotation.}
\end{lema}

\textbf{Proof.} Let $v_{n,z}$ and $v^{n,z}$ be frieze vectors, fixed a labeling $K_{1}$ in an $(n+3)$-polygon, then 

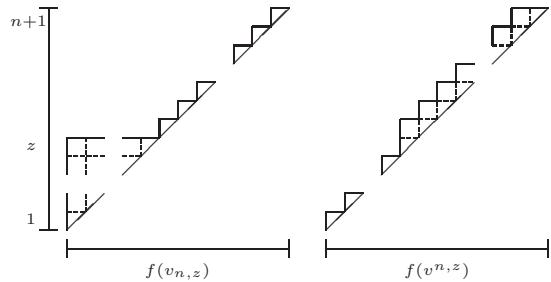
\begin{center}
	\begin{picture}(175, 144.5)
		
		\multiput(10.5,3.5)(98,0){2} {\line(1,0){84}}
		\multiput(10.5,0)(98,0){2} {\line(0,1){7}}
		\multiput(94.5,0)(98,0){2} {\line(0,1){7}}
		\multiput(0,10.5)(0,84){2} {\line(1,0){7}}
		\put(3.5,10.5){\line(0,1){84}}
		
		\multiput(10.5,10.5)(98,0){2}{\line(1,1){14}}
		\multiput(31.5,31.5)(98,0){2}{\line(1,1){35}}
		\multiput(73.5,73.5)(98,0){2}{\line(1,1){21}}
		
		\put(10.5, 10.5){\line(0,1){14}}
		\put(10.5, 31.5){\line(0,1){14}}
		\put(10.5, 45.5){\line(1,0){14}}
		\put(31.5, 45.5){\line(1,0){14}}
		\multiput(45.5,45.5)(28,28){2}{
			\multiput(0,0)(7,7){3}{\line(0,1){7}}
			\multiput(0,7)(7,7){3}{\line(1,0){7}}}
		\multiput(10.5,17.5)(2.5,0){3} {\line(1,0){1.5}}
		\multiput(10.5,38.5)(2.5,0){6} {\line(1,0){1.5}}
		\multiput(31.5,38.5)(2.5,0){3} {\line(1,0){1.5}}
		\multiput(17.5,17.5)(0,2.5){3} {\line(0,1){1.5}}
		\multiput(17.5,31.5)(0,2.5){6} {\line(0,1){1.5}}
		\multiput(38.5,38.5)(0,2.5){3} {\line(0,1){1.5}}
		
		\multiput(108.5,10.5)(7,7){2} {\line(0,1){7}}
		\multiput(108.5,17.5)(7,7){2} {\line(1,0){7}}
		\put(129.5,31.5){\line(0,1){7}}
		\put(129.5,38.5){\line(1,0){7}}
		\put(136.5,38.5){\line(0,1){14}}
		\multiput(136.5,52.5)(7,7){4} {\line(1,0){7}}
		\multiput(143.5,52.5)(7,7){3} {\line(0,1){7}}
		\multiput(171.5,80.5)(7,7){2} {\line(0,1){7}}
		\put(171.5,87.5){\line(1,0){7}}
		\put(178.5,94.5){\line(1,0){14}}
		\multiput(136.5,45.5)(7,7){3}{\multiput(0,0)(2.5,0){3} {\line(1,0){1.5}} \multiput(7,0)(0,2,5){3}{\line(0,1){1.5}}}
		\multiput(157.5,66.5)(2.5,0){3} {\line(1,0){1.5}}
		\multiput(171.5,80.5)(7,7){2}{\multiput(0,0)(2.5,0){3} {\line(1,0){1.5}} \multiput(7,0)(0,2,5){3}{\line(0,1){1.5}}}
		\put(-5,14){$_{\footnotesize{1}}$}
		\put(-5,42){$_{\footnotesize{z}}$}
		\put(-11,89){$_{\footnotesize{n+1}}$}
		
		\put(40,-5) {$_{f( v_{n,z})}$ }
		\put(138,-5) {$_{f(v^{n,z})}$ }
		
	\end{picture}
\end{center}

 the following identities hold by applying  map $F$ (see (\ref{F})) as follows:\\
 \begin{equation}
 \begin{split}
F(v_{n,z})&=(\underbrace{n}_{1},\dots,\underbrace{z}_{n-z},\underbrace{0}_{n-z+1},\dots,\underbrace{0}_{n}),\\ F(v^{n,z})&=(\underbrace{n-1}_{1},\dots,\underbrace{z-1}_{n-z},\underbrace{z-1}_{n-z+1},\dots,  \underbrace{1}_{n}), \notag
\end{split}
\end{equation}
if we change $K_{1}$  for $K_{2}$ by making the following replacements:
\begin{itemize}
\item Vertex $l \in K_{1}$ is changed for $l-1 \in K_{2}$, $1 \leq l \leq n+2$,
\item Vertex $0 \in K_{1}$ is changed for $n+2 \in K_{2}$,
\end{itemize}

then the diagonals from $0$ to $r_{1}$ in $K_{1}$ are diagonals from $r_{1}-1$ to $n+2$ in $K_{2}$,   and the diagonals from $r_{2}$ in $K_{1}$ are diagonals from $r_{2}-1$ in $K_{2}$, for $0 \leq r_{1} \leq z \leq r_{2}\leq n$. Therefore $ F(v_{n,z}) \in K_{1}$ coincides with $F(v^{n,z}) \in K_{2}$.$\hfill\square$\\

Note that, there exists a permutation

\footnotesize{\[ \sigma=
   \left ( 
      \begin{matrix} 
         1 & 2 & \dots & n-z-1 & n-z & n-z+1 & n-z+2 & \dots & n-1   & n \\
         1 & 2 & \dots & n-z-1 & n-z &  n    & n-1   & \dots & n-z+2 & n-z+1
         \end{matrix}
   \right)\]}
\small
in $S_{n}$ that describes a bijection between the coordinates of the vector $F(v_{n,z})=(u_{1}, \dots, u_{n})$ and the vector $F(v^{n,z})=(u'_{1}, \dots, u'_{n})$ such that $\sigma(F(v_{n,z}))=(u_{\sigma(1)}, \dots,$ $ u_{\sigma(n)})=(u'_{1}, \dots, u'_{n})=F(v^{n,z})$ in $K_{2}$. In general, if $v$ and $w$ realize the same  triangulation except for one anti-clockwise  rotation, then there exists a permutation $\sigma' \in S_{n}$ such that $\sigma'(F(v))=F(w)$ in $K_{2}$.

\begin{lema} \label{lema4.5}

\textit{ Let $A$ and $B$ be a coupling of positive integral diamonds of Dynkin type $\mathbb{A}_{n}$, and $v_{A}=(a_{1},\dots, a_{z}, \dots, a_{n})$ a corresponding associated vector with $a_{t}=1$ for $z\leq t\leq n$. If $v_{A}$ and $v_{B}$ realize the same triangulation except for one anti-clockwise rotation. Then vectors
\begin{equation}
\begin{split}
v_{A+i}&=(a_{1}, \dots, a_{i-1},a_{i-1}+a_{i},a_{i+1}, \dots, a_{n-1})\quad \text{and}\\
v_{B+i-1}&= (b_{1}, \dots , b_{i-2}, b_{i-2}+b_{i-1}, b_{i-1}, \dots, b_{n-1}),\notag
\end{split}
\end{equation}
realize the same triangulation except for one anti-clockwise rotation for $z-1 \leq i \leq n$, and $i \geq 2$.}
\end{lema} 
 
\textbf{Proof.} Let $v_{A}$ and $v_{B}$ be associated vectors to the $\mathbb{A}_{n}$-diamonds $A$ and $B$, respectively. Since $v_{A}$ and $v_{B}$  realize the same triangulation except for one anti-clockwise rotation then there exists a permutation $\sigma \in S_{n}$ such that $\sigma(F(v_{A}))=F(v_{B})$ in $K_{2}$. The following options arise from the map $f$, such that: 

\begin{enumerate}[(1)]
\item If $i> z \geq 1$, then
\begin{equation}
\begin{split}
 f(v_{A})&=(\dots, \underbrace{1}_{i-1},\underbrace{1}_{i},\dots ),\\
 f(v_{B})&=(\dots, \underbrace{d}_{i-2}, \underbrace{2}_{i-1}, \underbrace{2}_{i}, \dots),\quad\text{ (see Figure \ref{figure4.5})}\notag
\end{split}
\end{equation} 

\begin{enumerate}
\item [(1.1)] If $d=1$, then
\begin{equation}
\begin{split}
F(v_{A})&=(\dots,\underbrace{i}_{n-i},\underbrace{i-1}_{n+1-i},\dots ), \\\notag
F(v_{B})&=(\dots,\underbrace{i-1}_{n-i},\underbrace{i-2}_{n+1-i},\dots ),
\end{split}
\end{equation}
and $\sigma$ satisfies the expression,

\begin{equation} \label{equa4.1.4}
\sigma(r)=\begin{cases}
  r, & \mbox{if $ r \leq n+1-i$,}\\\notag
 m, & \mbox{otherwise,}
\end{cases}
\end{equation}
for some $m> n+1-i$. Applying  $F$ to $v_{A+i}$ and $v_{B+i-1}$, it holds that
\begin{equation}
\begin{split}
F(v_{A+i})&=(\dots,\underbrace{i-1}_{n-i},\underbrace{i-1}_{n+1-i},\dots ),\\ 
 F(v_{B+i-1})&=(\dots,\underbrace{i-2}_{n-i},\underbrace{i-2}_{n+1-i},\dots ),\notag
\end{split}
\end{equation} 
then there exits $\sigma' \in S_{n}$ such that  $\sigma'=\sigma$ and $\sigma'(F(v_{A+i}))=F(v_{B+i-1})$ in $K_{2}$ (see Figure \ref{figure4.5}).\\

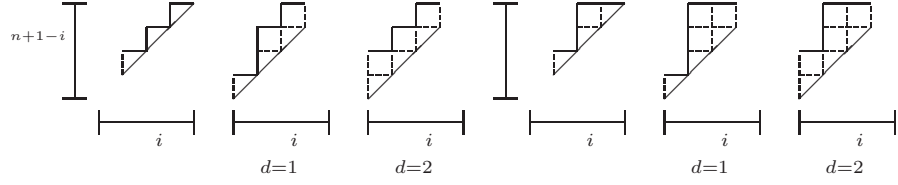
\begin{figure}[h!]
	\begin{center}
		\begin{picture}(290, 50)
			
			\put(-20,38){$_{_{n+1-i}}$}
			
			\multiput(0,0)(163,0){2}{
				\put(4.5,13.5){\line(0,1){36}}
				\multiput(13.5,0)(51,0){3}{\line(0,1){9}}
				\multiput(49.5,0)(51,0){3}{\line(0,1){9}}
				\multiput(13.5,4.5)(51,0){3} {\line(1,0){36}}
				\multiput(0,13.5)(0,36){2}{\line(1,0){9}}
				\multiput(35,-3)(51,0){3}{$_{i}$}
				\put(22.5,22.5){\line(1,1){27}}
				\multiput(64.5,13.5)(51,0){2}{\line(1,1){27}}
				\put(74.5,-13){$_{d=1}$}
				\put(125.5,-13){$_{d=2}$}
			}
			
			\multiput(22.5,31.5)(9,9){3}{\line(1,0){9}}
			\multiput(31.5,31.5)(9,9){2}{\line(0,1){9}}
			\multiput(22.5,22.5)(0,2.9){3} {\line(0,1){1.9}}
			
			\put(64.5,22.5){\line(1,0){9}}
			\put(73.5,22.5){\line(0,1){18}}
			\multiput(73.5,40.5)(9,9){2}{\line(1,0){9}}
			\put(82.5,40.5){\line(0,1){9}}
			\multiput(64.5,13.5)(0,2.9){3} {\line(0,1){1.9}}
			\multiput(73.5,31.5)(9,9){2}{\multiput(0,0)(2.9,0){3}{\line(1,0){1.9}} \multiput(9,0)(0,2.9){3}{\line(0,1){1.9}}}
			
			\multiput(115.5,31.5)(9,9){3}{\line(1,0){9}}
			\multiput(124.5,31.5)(9,9){2}{\line(0,1){9}}
			\multiput(115.5,13.5)(0,2.9){6}{\line(0,1){1.9}}
			\multiput(115.5,22.5)(9,9){3}{\multiput(0,0)(2.9,0){3}{\line(1,0){1.9}} \multiput(9,0)(0,2.8){3}{\line(0,1){1.8}}}
			
			\multiput(185.5,31.5)(93,0){2}{\put(0,0){\line(1,0){9}}
				\put(9,0){\line(0,1){18}}
				\put(9,18){\line(1,0){18}}}
			
			\multiput(185.5,22.5)(18,18){2}{\multiput(0,0)(0,2.9){3}{\line(0,1){1.8}}}
			\multiput(194.5,40.5)(2.9,0){3}{\line(1,0){1.9}}
			
			\put(227.5,22.5){\line(1,0){9}}
			\put(236.5,22.5){\line(0,1){27}}
			\put(236.5,49.5){\line(1,0){18}}
			\multiput(227.5,13.5)(27,27){2}{\multiput(0,0)(0,2.9){3}{\line(0,1){1.9}}}
			\multiput(245.5,31.5)(0,2.9){6}{\line(0,1){1.9}}
			\multiput(236.5,31.5)(2.9,0){3}{\line(1,0){1.9}}
			\multiput(236.5,40.5)(2.9,0){6}{\line(1,0){1.9}}
			
			\multiput(278.5,13.5)(18,18){2}{\multiput(0,0)(0,2.9){6}{\line(0,1){1.9}}}
			\multiput(287.5,22.5)(18,18){2}{\multiput(0,0)(0,2.9){3}{\line(0,1){1.9}}}
			\multiput(278.5,22.5)(9,9){2}{\multiput(0,0)(2.9,0){3}{\line(1,0){1.9}}}
			\multiput(287.5,40.5)(2.9,0){6}{\line(1,0){1.9}}
		\end{picture}
		\par\bigskip
	\end{center}
	\caption{Dyck paths associated to vectors $v_{A}$, $v_{B}$ (left), $v_{A+i}$ and $v_{B+i-1}$ (right)  for $i>z$.}\label{figure4.5}
\end{figure}

\item [(1.2)] The case for $d=2$ is the same as the previous case.
\item [(1.3)]If $d=3$, then $A$ and $B$  do not realize the same triangulation.
\end{enumerate}

Note that, if $z=1$, this case satisfies the condition (1.1) and (1.2) without $d$.

\item If $i= z \geq 2$,  then
\begin{equation}
\begin{split}
f(v_{A})&=(\dots, \underbrace{2}_{i-1},\underbrace{1}_{i},\dots ),\\\notag
f(v_{B})&=(\dots, \underbrace{b}_{i-2}, \underbrace{a}_{i-1}, \underbrace{2}_{i}, \dots), \quad\text{ (see Figure \ref{figure4.7})}.
\end{split}
\end{equation}
\begin{enumerate}
\item [(2.1)] If $a=1$ and $b=1$, then 

\begin{equation}
\begin{split}
F(v_{A})&=(\dots,\underbrace{i}_{n-i},\dots ),\\\notag
F(v_{B})&=(\dots,\underbrace{i-1}_{n-i},\underbrace{i-1}_{n+1-i},\dots ),
\end{split}
\end{equation}
and $\sigma_{1}$ is defined by the following cases:

\begin{equation}
\sigma_{1}(r)=\begin{cases}
  r, & \mbox{if $ r \leq n-i$,}\\\notag
  n+1-i, & \mbox{if $  r =n$,}\\
 m, & \mbox{otherwise,}
\end{cases}
\end{equation}

for some $m> n+1-i$.  Applying $F$, we obtain 
\begin{equation}
\begin{split}
F(v_{A+i})&=(\dots,\underbrace{i-1}_{n-i},\dots ),\\ 
F(v_{B+i-1})&=(\dots,\underbrace{i}_{n-i},\underbrace{i-2}_{n+1-i},\dots ),\notag
\end{split}
\end{equation}
 then there exits $\sigma'_{1} \in S_{n}$ satisfying the following cases:  

\begin{equation}
\sigma'_{1}(r)=\begin{cases}
   n-i, & \mbox{if $ r=n$,}\\
  n+1-i, & \mbox{if $  r =n-i$,}\\\notag
 \sigma_{1}(r), & \mbox{otherwise,}
\end{cases}
\end{equation}

 therefore $\sigma'_{1}(F(v_{A+i}))=F(v_{B+i-1})$ in $K_{2}$ (see Figure \ref{figure4.8}).\\

\begin{figure}[h!]
	\begin{center}
		\begin{picture}(180,50)
			
			\put(-20,38){$_{_{n+1-i}}$}
			\put(4.5,13.5){\line(0,1){36}}
			\multiput(13.5,0)(51,0){4}{\line(0,1){9}}
			\multiput(49.5,0)(51,0){4}{\line(0,1){9}}
			\multiput(13.5,4.5)(51,0){4} {\line(1,0){36}}
			\multiput(0,13.5)(0,36){2}{\line(1,0){9}}
			\multiput(35,-3)(51,0){4}{$_{i}$}
			\put(22.5,22.5){\line(1,1){27}}
			\multiput(64.5,13.5)(51,0){3}{\line(1,1){27}}
			\put(74.5,-13){$_{a=1}$}
			\put(125.5,-13){$_{a=1}$}
			\put(176.5,-13){$_{a=2}$}
			\put(74.5,-21){$_{b=1}$}
			\put(125.5,-21){$_{b=2}$}
			\put(176.5,-21){$_{b=3}$}

			\put(22.5,40.5){\line(1,0){18}}
			\put(40.5,40.5){\line(0,1){9}}
			\put(40.5,49.5){\line(1,0){9}}
			\multiput(22.5,22.5)(0,2.9){6}{\line(0,1){1.9}}
			\multiput(31.5,31.5)(0,2.9){3}{\line(0,1){1.9}}
			\multiput(22.5,31.5)(2.9,0){3}{\line(1,0){1.9}}
			
			\multiput(64.5,22.5)(9,9){2}{\put(0,0){\line(1,0){9}}}
			\put(73.5,22.5){\line(0,1){9}}
			\put(82.5,31.5){\line(0,1){18}}
			\put(82.5,49.5){\line(1,0){9}}
			\multiput(64.5,13.5)(27,27){2}{\multiput(0,0)(0,2.9){3}{\line(0,1){1.9}}}
			\multiput(82.5,41.5)(2.9,0){3}{\line(1,0){1.9}}
			
			\put(115.5,31.5){\line(1,0){18}}
			\put(133.5,31.5){\line(0,1){18}}
			\put(133.5,49.5){\line(1,0){9}}
			\multiput(115.5,22.5)(18,18){2}{\multiput(0,0)(2.9,0){3}{\line(1,0){1.9}}
				\multiput(9,0)(0,2.9){3}{\line(0,1){1.9}}}
			\multiput(115.5,13.5)(0,2.9){6}{\line(0,1){1.9}}
			
			\put(166.5,40.5){\line(1,0){18}}
			\put(184.5,40.5){\line(0,1){9}}
			\put(184.5,49.5){\line(1,0){9}}
			\multiput(166.5,13.5)(0,2.9){9}{\line(0,1){1.9}}
			\multiput(175.5,22.5)(0,2.9){6}{\line(0,1){1.9}}
			\multiput(184.5,31.5)(9,9){2}{\multiput(0,0)(0, 2.9){3}{\line(0,1){1.9}}}
			\multiput(166.5,22.5)(18,18){2}{\multiput(0,0)(2.9,0){3}{\line(1,0){1.9}}}
			\multiput(166.5,31.5)(2.9,0){6}{\line(1,0){1.9}}
		\end{picture}
	\end{center}
	\par\bigskip
	\par\bigskip
	\caption{Dyck paths associated to vectors $v_{A}$ and $v_{B}$  for $i=z$.}\label{figure4.7}
\end{figure}
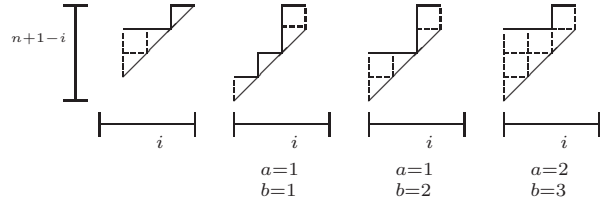

\item [(2.2)] If $a=1$ and $b=2$, then conditions defined in the case (2.1) hold.
\item [(2.3)] For $a=2$ and $b=1$ or $b=2$, we have only contradictions.
\item [(2.4)] If $a=2$ and $b=3$, $F(v_{B})=(\dots,\underbrace{i-1}_{n-i},\dots )$ and $\sigma_{2}=\sigma$. Applying $F$ to $v_{B+i-1}$, it holds that $ F(v_{B+i-1})=(\dots,\underbrace{i-2}_{n-i},\dots )$ then there exits $\sigma'_{2} \in S_{n}$ such that  $\sigma'_{2}=\sigma$ and $\sigma'_{2}(F(v_{A+i}))=F(v_{B+i-1})$ in $K_{2}$ (see Figure \ref{figure4.8}).
\item [(2.5)] Case (2.3) holds for $a=3$ and $b=1,2,3$.

\begin{figure}[h!]
	
	\begin{center}
		\begin{center}
			\begin{picture}(180,50)
				
				\put(-20,38){$_{_{n+1-i}}$}
				\put(4.5,13.5){\line(0,1){36}}
				\multiput(13.5,0)(51,0){4}{\line(0,1){9}}
				\multiput(49.5,0)(51,0){4}{\line(0,1){9}}
				\multiput(13.5,4.5)(51,0){4} {\line(1,0){36}}
				\multiput(0,13.5)(0,36){2}{\line(1,0){9}}
				\multiput(35,-3)(51,0){4}{$_{i}$}
				\put(22.5,22.5){\line(1,1){27}}
				\multiput(64.5,13.5)(51,0){3}{\line(1,1){27}}
				\put(74.5,-13){$_{a=1}$}
				\put(125.5,-13){$_{a=1}$}
				\put(176.5,-13){$_{a=2}$}
				\put(74.5,-21){$_{b=1}$}
				\put(125.5,-21){$_{b=2}$}
				\put(176.5,-21){$_{b=3}$}
				
				\put(22.5,40.5){\line(1,0){9}}
				\put(31.5,40.5){\line(0,1){9}}
				\put(31.5,49.5){\line(1,0){18}}
				\multiput(22.5,31.5)(9,9){2}{\multiput(9,0)(0,2.9){3}{\line(0,1){1.9}}
					\multiput(0,0)(2.9,0){3}{\line(1,0){1.9}}}
				\multiput(22.5,22.5)(0,2.8){6}{\line(0,1){1.8}}
				
				\put(64.5,22.5){\line(1,0){9}}
				\put(73.5,22.5){\line(0,1){18}}
				\put(73.5,40.5){\line(1,0){18}}
				\put(91.5,40.5){\line(0,1){9}}
				\multiput(64.5,13.5)(18,18){2}{\multiput(0,0)(0,2.9){3}{\line(0,1){1.9}}}
				\multiput(73.5,31.5)(2.9,0){3}{\line(1,0){1.9}}
				
				\put(115.5,31.5){\line(1,0){9}}
				\multiput(124.5,31.5)(18,9){2}{\put(0,0){\line(0,1){9}}}
				\put(124.5,40.5){\line(1,0){18}}
				\multiput(115.5,22.5)(9,9){2}{\multiput(0,0)(2.9,0){3}{\line(1,0){1.9}} \multiput(9,0)(0,2.9){3}{\line(0,1){1.9}}}
				\multiput(115.5,13.5)(0,2.9){6}{\line(0,1){1.9}}
				
				\put(166.5,40.5){\line(1,0){9}}
				\put(175.5,40.5){\line(0,1){9}}
				\put(175.5,49.5){\line(1,0){18}}
				\multiput(166.5,13.5)(0,2.9){9}{\line(0,1){1.9}}
				\multiput(175.5,22.5)(9,9){2}{\multiput(0,0)(0,2.9){6}{\line(0,1){1.9}}}
				\multiput(193.5,40.9)(0,2.9){3}{\line(0,1){1.9}}
				\multiput(166.5,31.5)(9,9){2}{\multiput(0,0)(2.9,0){6}{\line(1,0){1.9}}}
				\multiput(166.5,22.5)(2.9,0){3}{\line(1,0){1.9}}
			\end{picture}
		\end{center}
		\par\bigskip
		\par\bigskip
	\end{center}
	\caption{Dyck paths associated to vectors $v_{A+i}$ and $v_{B+i-1}$ for $i=z$.}\label{figure4.8}
\end{figure}
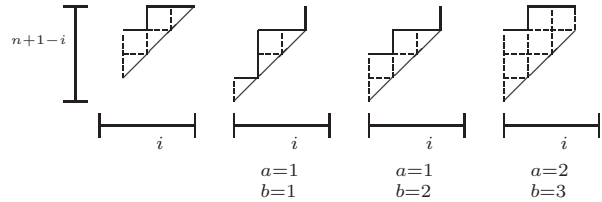
\end{enumerate}

Note that, if $z=2$ then conditions for $a=1,2$ without $b$ hold.

\item [(3)] If $i=z-1 \geq 3$, then
\begin{equation}
\begin{split}
f(v_{A})&=(\dots, \underbrace{2}_{i},\underbrace{1}_{i+1},\dots ),\\\notag
f(v_{B})&=(\dots, \underbrace{b}_{i-1}, \underbrace{a}_{i}, \underbrace{2}_{i+1}, \dots)\quad\text{ (see  Figure \ref{figure4.9})}.
\end{split}
\end{equation}
\begin{enumerate}
\item [(3.1)] If $a=1$ and $b=1$, then
\begin{equation}
\begin{split}
F(v_{A})&=(\dots,\underbrace{i+1}_{n-i-1},\dots ),\\\notag 
F(v_{B})&=(\dots,\underbrace{i}_{n-i-1},\underbrace{i}_{n-i}, \underbrace{i-1}_{n+1-i},\dots ),\quad\text{and}\\  
\end{split}
\end{equation}
\begin{equation}
\sigma_{3}(r)=\begin{cases}
  r, & \mbox{if $ r \leq n-i-1$,}\\
  n-i, & \mbox{if $  r =n$,}\\
  n+1-i, & \mbox{if $  r =n-1$,}\\\notag
 m, & \mbox{otherwise,}
\end{cases}
\end{equation}

for some $m> n+1-i$. Provided that 
\begin{equation}
\begin{split}
F(v_{A+i})&=(\dots,\underbrace{i-1}_{n-i-1},\dots ),\quad \text{ and }\\
F(v_{B+i-1})&=(\dots,\underbrace{i+1}_{n-i-1},\underbrace{i}_{n-i},\underbrace{i-2}_{n+1-i},\dots ),\notag
\end{split}
\end{equation}
 then, there exits $\sigma'_{3} \in S_{n}$ such that 

\begin{equation}
\sigma'_{3}(r)=\begin{cases}
   n-i-1, & \mbox{if $ r=n$,}\\
     n-i, & \mbox{if $ r=n-1$,}\\
  n+1-i, & \mbox{if $  r =n-i-1$,}\\\notag
 \sigma_{3}(r), & \mbox{otherwise,}
\end{cases}
\end{equation}

 then $\sigma'_{3}(F(v_{A+i}))=F(v_{B+i-1})$ in $K_{2}$ (see Figure \ref{figure4.10}).\\
 
\begin{figure}[h!]
	\begin{center}
		
		\begin{picture}(180,50)
			\put(-20,29){$_{_{n+1-i}}$}
			\put(4.5,13.5){\line(0,1){36}}
			\multiput(13.5,0)(51,0){4}{\line(0,1){9}}
			\multiput(49.5,0)(51,0){4}{\line(0,1){9}}
			\multiput(13.5,4.5)(51,0){4} {\line(1,0){36}}
			\multiput(0,13.5)(0,36){2}{\line(1,0){9}}
			\multiput(26,-3)(51,0){4}{$_{i}$}
			\put(22.5,22.5){\line(1,1){27}}
			\multiput(64.5,13.5)(51,0){3}{\line(1,1){27}}
			\put(74.5,-13){$_{a=1}$}
			\put(125.5,-13){$_{a=1}$}
			\put(176.5,-13){$_{a=2}$}
			\put(74.5,-21){$_{b=1}$}
			\put(125.5,-21){$_{b=2}$}
			\put(176.5,-21){$_{b=3}$}

			\put(22.5,40.5){\line(1,0){18}}
			\put(40.5,40.5){\line(0,1){9}}
			\put(40.5,49.5){\line(1,0){9}}
			\multiput(22.5,22.5)(0,2.9){6}{\line(0,1){1.9}}
			\multiput(31.5,31.5)(0,2.9){3}{\line(0,1){1.9}}
			\multiput(22.5,31.5)(2.9,0){3}{\line(1,0){1.9}}
			
			\multiput(64.5,22.5)(9,9){2}{\put(0,0){\line(1,0){9}}}
			\put(73.5,22.5){\line(0,1){9}}
			\put(82.5,31.5){\line(0,1){18}}
			\put(82.5,49.5){\line(1,0){9}}
			\multiput(64.5,13.5)(27,27){2}{\multiput(0,0)(0,2.9){3}{\line(0,1){1.9}}}
			\multiput(82.5,41.5)(2.9,0){3}{\line(1,0){1.9}}
			
			\put(115.5,31.5){\line(1,0){18}}
			\put(133.5,31.5){\line(0,1){18}}
			\put(133.5,49.5){\line(1,0){9}}
			\multiput(115.5,22.5)(18,18){2}{\multiput(0,0)(2.9,0){3}{\line(1,0){1.9}}
				\multiput(9,0)(0,2.9){3}{\line(0,1){1.9}}}
			\multiput(115.5,13.5)(0,2.9){6}{\line(0,1){1.9}}
			
			\put(166.5,40.5){\line(1,0){18}}
			\put(184.5,40.5){\line(0,1){9}}
			\put(184.5,49.5){\line(1,0){9}}
			\multiput(166.5,13.5)(0,2.9){9}{\line(0,1){1.9}}
			\multiput(175.5,22.5)(0,2.9){6}{\line(0,1){1.9}}
			\multiput(184.5,31.5)(9,9){2}{\multiput(0,0)(0, 2.9){3}{\line(0,1){1.9}}}
			\multiput(166.5,22.5)(18,18){2}{\multiput(0,0)(2.9,0){3}{\line(1,0){1.9}}}
			\multiput(166.5,31.5)(2.9,0){6}{\line(1,0){1.9}}
		\end{picture}
	\end{center}
	\par\bigskip
	\par\bigskip
	\caption{Dyck paths associated to vectors $v_{A}$ and $v_{B}$  for $i-1=z$.}\label{figure4.9}
\end{figure}
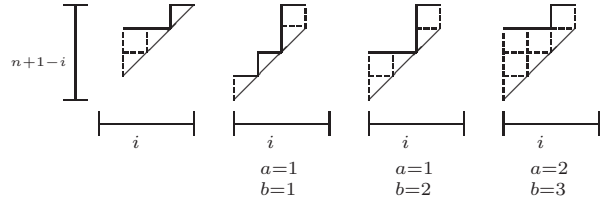 
 
\item [(3.2)] If $a=1$ and $b=2$. It holds that, 
\begin{equation}
\begin{split}
F(v_{B})&=(\dots,\underbrace{i}_{n-i-1},\underbrace{i}_{n-i},\dots ),\notag
\end{split}
\end{equation}

whereas, $\sigma_{4}$ is given  by the identities

\begin{equation}
\sigma_{4}(r)=\begin{cases}
  r, & \mbox{if $ r \leq n-i-1$,}\\
  n-i, & \mbox{if $  r =n$,}\\\notag
 m, & \mbox{otherwise,}
\end{cases}
\end{equation}
for some $m>n-i$. Applying $F$ to $v_{B+i-1}$, we get
\begin{equation}
\begin{split}
F(v_{B+i-1})&=(\dots,\underbrace{i+1}_{n-i-1},\underbrace{i-2}_{n-i},\dots ),\notag
\end{split}
\end{equation}
 then there exists $\sigma'_{4}$ with

\begin{equation}
\sigma'_{4}(r)=\begin{cases}
  n-i-1, & \mbox{if $ r \leq n$,}\\\notag
  n-i, & \mbox{if $  r =n-i-1$,}\\
 \sigma_{4}(r), & \mbox{otherwise,}
\end{cases}
\end{equation}

therefore $\sigma'_{4}(F(v_{A+i}))=F(v_{B+i-1})$ in $K_{2}$ (see Figure \ref{figure4.10}).
\item [(3.3)] If $a=2$ and $b=3$. $F(v_{B})= (\dots,\underbrace{i}_{n-i-1},\dots )$, provided that 

\begin{equation}
\sigma_{5}(r)=\begin{cases}
  r, & \mbox{if $ r \leq n-i-1$,}\\\notag
 m, & \mbox{otherwise,}
\end{cases}
\end{equation}
for some $m>n-i$. In this case, $F(v_{B+i-1})=(\dots,\underbrace{i-2}_{n-i-1},\dots )$, and there is $\sigma'_{5}=\sigma_{5}$ such that $\sigma'_{5}(F(v_{A+i}))=F(v_{B+i-1})$ in $K_{2}$ (see Figure \ref{figure4.10}).

\end{enumerate}

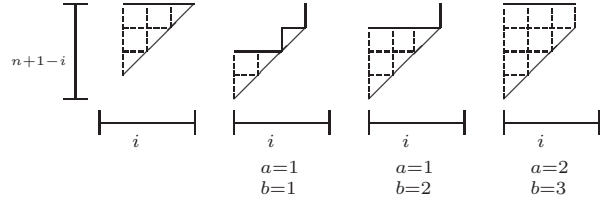
\begin{figure}[h!]
	
	\begin{center}
		\begin{picture}(180,50)
			\put(-20,29){$_{_{n+1-i}}$}
			\put(4.5,13.5){\line(0,1){36}}
			\multiput(13.5,0)(51,0){4}{\line(0,1){9}}
			\multiput(49.5,0)(51,0){4}{\line(0,1){9}}
			\multiput(13.5,4.5)(51,0){4} {\line(1,0){36}}
			\multiput(0,13.5)(0,36){2}{\line(1,0){9}}
			\multiput(26,-3)(51,0){4}{$_{i}$}
			\put(22.5,22.5){\line(1,1){27}}
			\multiput(64.5,13.5)(51,0){3}{\line(1,1){27}}
			\put(74.5,-13){$_{a=1}$}
			\put(125.5,-13){$_{a=1}$}
			\put(176.5,-13){$_{a=2}$}
			\put(74.5,-21){$_{b=1}$}
			\put(125.5,-21){$_{b=2}$}
			\put(176.5,-21){$_{b=3}$}
			
			\put(22.5,49.5){\line(1,0){27}}
			\multiput(22.5,22.5)(0,2.9){9}{\line(0,1){1.9}}
			\multiput(31.5,31.5)(0,2.9){6}{\line(0,1){1.9}}
			\multiput(40.5,40.5)(0,2.9){3}{\line(0,1){1.9}}
			\multiput(22.5,31.5)(2.9,0){3}{\line(1,0){1.9}}
			\multiput(22.5,40.5)(2.9,0){6}{\line(1,0){1.9}}
			
			\put(64.5,31.5){\line(1,0){18}}
			\multiput(82.5,31.5)(9,9){2}{\put(0,0){\line(0,1){9}}}
			\put(82.5,40.5){\line(1,0){9}}
			\multiput(64.5,13.5)(0,2.9){6}{\line(0,1){1.9}}
			\multiput(64.5,22.5)(2.9,0){3}{\line(1,0){1.9}}
			\multiput(73.5,22.5)(0,2.9){3}{\line(0,1){1.9}}
			
			\put(115.5,40.5){\line(1,0){27}}
			\put(142.5,40.5){\line(0,1){9}}
			\multiput(115.5,13.5)(0,2.9){9}{\line(0,1){1.9}}
			\multiput(124.5,22.5)(0,2.9){6}{\line(0,1){1.9}}
			\multiput(132.5,31.5)(0,2.9){3}{\line(0,1){1.9}}
			\multiput(115.5,22.5)(2.9,0){3}{\line(1,0){1.9}}
			\multiput(115.5,31.5)(2.9,0){6}{\line(1,0){1.9}}
			
			\put(166.5,49.5){\line(1,0){27}}
			\multiput(166.5,13.5)(0,2.9){12}{\line(0,1){1.9}}
			\multiput(175.5,22.5)(0,2.9){9}{\line(0,1){1.9}}
			\multiput(184.5,31.5)(0,2.9){6}{\line(0,1){1.9}}
			\multiput(193.5,40.5)(0,2.9){3}{\line(0,1){1.9}}
			\multiput(166.5,22.5)(2.9,0){3}{\line(1,0){1.9}}
			\multiput(166.5,31.5)(2.9,0){6}{\line(1,0){1.9}}
			\multiput(166.5,40.5)(2.9,0){9}{\line(1,0){1.9}}
		\end{picture}
	\end{center}
	\par\bigskip
	\par\bigskip
	\caption{Dyck paths associated to vectors $v_{A+i}$ and $v_{B+i-1}$  for $i-1=z$.}\label{figure4.10}
\end{figure} 

Same arguments are used for the remaining cases (see item (2) of this proof).$\hfill\square$
\end{enumerate}

\addtocounter{prop}{5}
\begin{prop} \label{prop4.9}
 \textit{Two positive integral diamonds  of Dynkin type $\mathbb{A}_{n}$ are in the same minimal $p$-cycle   if their triangulations  are in the same mutation  class.}
\end{prop}
 
\textbf{Proof.} It is a direct consequence of Proposition \ref{prop4.3}, Lemmas \ref{lema4.4} and \ref{lema4.5}$\hfill\square$\\
 
The following result gives a way to build frieze patterns. 

\addtocounter{teor}{3}
\begin{teor}\label{main}
\textit{Let $A^{0}$ be a positive integral diamond of Dynkin type $\mathbb{A}_{n}$ and let  $\lbrace A^{t} \rbrace_{0 \leq t \leq p-1}$ be the minimal $p$-cycle  generated by $A^{0}$. Then:}

\begin{enumerate}
\item  [(i) ]\textit{$A^{0}$ and $F(v_{A^{0}})$ generate the same frieze pattern (see $\mathrm{(\ref{F})}$).}
\item  [(ii)]\textit{$\lbrace A^{t}\rbrace_{0 \leq t \leq p-1}$ is in  surjective correspondence with a direct sum of $p$ indecomposable objects of a Dyck paths category.}
\end{enumerate}

\end{teor}
 
\textbf{Proof.} Let $\mathbb{D}_{\mathbb{A}_{n}}$ be the set of all vectors associated to positive integral diamonds of Dynkin type $\mathbb{A}_{n}$, $A^{0}$ a positive integral diamond of Dynkin type  $\mathbb{A}_{n}$, and $\lbrace A^{t} \rbrace_{0 \leq t \leq p-1}$  the minimal $p$-cycle  generated by $A^{0}$.

\begin{enumerate}
\item [(i)] Let $K$ be a labeling of an $(n+3)$-polygon, Theorem \ref{teor4.2} implies that 

\begin{equation}
\begin{array}{rcl}
F(v_{A^{0}})&=& g((a^{0}_{11}, T_{2}(v_{A^{0}}), \dots, T_{n}(v_{A^{0}})))\\
&=& g(\lambda_{(a^{0}_{11}, T_{2}(v_{A^{0}}), \dots, T_{n}(v_{A^{0}}))})\\
&=& g((\lambda_{1}, \dots, \lambda_{n+1-a^{0}_{11}}, \underbrace{0, \dots, 0}_{a^{0}_{11}}))\\
&=& (l_{1}^{v_{1}}, \dots, l_{n+1-a_{11}^{0}}^{v_{n+1-a^{0}_{11}}}, l^{0}_{n-a^{0}_{11}}, \dots, l^{0}_{n}), \notag
\end{array}
\end{equation}
then, there are $a^{0}_{11}-1$ diagonals from the vertex $0$ to other vertices, i.e., there are $a^{0}_{11}$ triangles incident with vertex 0. Proposition \ref{prop4.9} allows us to establish that $a^{i}_{11}$ is the number of triangles incident with the vertex $i$, for $1 \leq i \leq n+3$, $i=pm$ and $1 \leq m \leq p \text{ }| \text{ }(n+3)$. Therefore $A^{0}$ and $F(v_{A^{0}})$ generate the same frieze pattern.

\item [(ii)] Let $(\mathfrak{D}_{2(n+1)},R)$ be  any Dyck paths category, we take objects of $(\mathfrak{D}_{2(n+1)},R)$ defined by the following identity

\begin{equation}
\overline{Ob} \text{ }(\mathfrak{D}_{2n},R)=\bigg \{ \bigoplus_{ G_{i} \in \mathfrak{D}_{2n}} G_{i} \text{ } \bigg | \text{ } g(\lambda_{G_{i}}) \text{ and } g(\lambda_{G_{j}}) \text{ are in the same mutation class} \bigg \},\notag
\end{equation}

 the map $\varphi: \mathbb{D}_{\mathbb{A}_{n}} \rightarrow  \overline{Ob} \text{ }(\mathfrak{D}_{2n},R)$, such that \[\varphi (v_{A^{0}})= f(v_{A^{0}})\oplus \dots \oplus f(v_{A^{p-1}}),\] with $\lbrace A^{t}\rbrace_{0 \leq t \leq p-1}$ is surjective as a consequence of Theorem \ref{teor4.2} and  Proposition \ref{prop4.9}. $\hfill\square$\\
\end{enumerate}

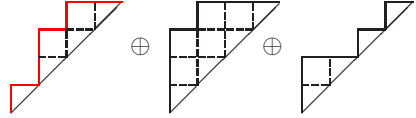
\begin{figure}[h!]
	\begin{center}
		\begin{picture}(94.5,44.5)
			
			\put(-30,0){\line(1,1){42}}
			\color{red}\put(-30,0){\line(0,1){10.5}}
			\put(-30,10.5){\line(1,0){10.5}}
			\put(-19.5,10.5){\line(0,1){21}}
			\put(-19.5,31.5){\line(1,0){10.5}}
			\put(-9,31.5){\line(0,1){10.5}}
			\put(-9,42){\line(1,0){21}}

			\normalcolor
			\multiput(-19.5,21)(3.6,0){3} {\line(1,0){2.6}} 
			\multiput(-9,31.5)(3.6,0){3} {\line(1,0){2.6}}
			\multiput(-9,21)(0,3.6){3} {\line(0,1){2.6}}
			\multiput(1.5,31.5)(0,3.6){3} {\line(0,1){2.6}}

			\put(30,0){\line(1,1){41.5}}
			\put(30,0){\line(0,1){31.5}}
			\put(30,31.5){\line(1,0){10.5}}
			\put(40.5,31.5){\line(0,1){10.5}}
			\put(40.5,42){\line(1,0){31}}
			
			\multiput(30,10.5)(3.6,0){3} {\line(1,0){2.6}} 
			\multiput(30,21)(3.6,0){6} {\line(1,0){2.6}} 
			\multiput(30,31.5)(3.6,0){9} {\line(1,0){2.6}}
			\multiput(40.5,10.5)(0,3.6){6} {\line(0,1){2.6}}
			\multiput(51,21)(0,3.6){6} {\line(0,1){2.6}}
			\multiput(61.5,31.5)(0,3.6){3} {\line(0,1){2.6}}

			
			\put(80,0){\line(1,1){42}}
			\put(80,0){\line(0,1){21}}
			\put(80,21){\line(1,0){21}}
			\put(101,21){\line(0,1){10.5}}
			\put(101,31.5){\line(1,0){10.5}}
			\put(111.5,31.5){\line(0,1){10.5}}
			\put(111.5,42){\line(1,0){10.5}}
			
			\multiput(80,10.5)(3.6,0){3} {\line(1,0){2.6}} 
			\multiput(90.5,10.5)(0,3.6){3} {\line(0,1){2.6}}

			\put(15,22.5){\textbf{$\oplus$}}
			\put(65,22.5){$\oplus$}
		\end{picture}
	\end{center}
	\caption{Examples of objects of a Dyck paths category.}\label{figure4.11}
\end{figure}

As an example of the use of Theorem \ref{main}, we choose the object $D$ of a Dyck paths category $(\mathfrak{D}_{2(n+1)},R)$ shown in Figure \ref{figure4.11}. Then $D$ has associated  the following frieze pattern where red numbers is a positive integral diamond of Dynkin type $\mathbb{A}_{n}$.

\[\xymatrix@=1mm{ &  \dots &&\textcolor{red}{1}&& 1 &&1 &&\textcolor{red}{1}&&1 &&\dots &\\
\dots&&\textcolor{red}{1}&&\textcolor{red}{3} && 2 && \textcolor{red}{1} && \textcolor{red}{3} && 2 && \dots  \\
&  \dots && \textcolor{red}{2} && \textcolor{red}{5} &&1 &&\textcolor{red}{2}&&\textcolor{red}{5} &&\dots &\\
\dots&& 1 &&\textcolor{red}{3} && \textcolor{red}{2} && 1 && \textcolor{red}{3} && \textcolor{red}{2} && \dots  \\
&  \dots && 1 &&  \textcolor{red}{1}&&1 && 1&&  \textcolor{red}{1}&&\dots &\\
}\]


\par\bigskip
Agust\'{\i}n Moreno Ca\~{n}adas\\
amorenoca@unal.edu.co\\
Department of Mathematics\\
Universidad Nacional de Colombia\\
Bogot\'a-Colombia\\

Isa\'{\i}as David Mar\'{\i}n  Gaviria\\
imaringa@unal.edu.co\\
Department of Mathematics\\
Universidad Nacional de Colombia\\
Bogot\'a-Colombia\\

Gabriel Bravo Rios\\
gbravor@unal.edu.co\\
Department of Mathematics\\
Universidad Nacional de Colombia\\
Bogot\'a-Colombia\\

Pedro Fernando Fern\'{a}ndez Espinosa\\
pffernandeze@unal.edu.co\\
Department of Mathematics\\
Universidad Nacional de Colombia\\
Bogot\'a-Colombia\\

\end{document}